\documentclass[12pt]{amsart}
\usepackage{amsmath,amsfonts,amssymb,amsthm}
\usepackage{graphicx,color}
\DeclareMathRadical{\sqrtsign}{symbols}{"70}{largesymbols}{"70}
\newcommand{\bb}{\mathbb}

\newcommand{\natls}{{\bb N}}

\newcommand{\reals}{{\bb R}}

\newlength{\figboxwidth}
\renewcommand{\bold}[1]{\medskip \noindent {\bf #1 }\nopagebreak}

\newcommand{\cross}{\times}

\newcommand{\st}{\;\: : \;\:}         

\def\@ifundefined#1#2#3%
  {\expandafter\ifx\csname#1\endcsname\relax#2\else#3\fi}

\@ifundefined{theoremstyle}{
}{
\theoremstyle{plain} 
}
\newtheorem{theorem}{Theorem}[section]

\newtheorem{proposition}[theorem]{Proposition}
\newtheorem{lemma}[theorem]{Lemma}

\@ifundefined{theoremstyle}{
}{
\theoremstyle{definition} 
}
\newtheorem{definition}[theorem]{Definition}

\newtheorem{remark}[theorem]{Remark}

\newcommand{\cK}{{\mathcal K}}

\newcommand{\cM}{{\mathcal M}}
\newcommand{\cN}{{\mathcal N}}

\newcommand{\cS}{{\mathcal S}}

\newcommand{\cV}{{\mathcal V}}

\mathchardef\GG="321D
\newcommand{\bA}{{\mathbf A}}
\newcommand{\bB}{{\mathbf B}}

\newcommand{\bH}{{\mathbf H}}

\newcommand{\bK}{{\mathbf K}}

\newcommand{\bM}{{\mathbf M}}
\newcommand{\bN}{{\mathbf N}}

\newcommand{\bP}{{\mathbf P}}

\newcommand{\Nbhd}{\operatorname{Nbhd}}

\newcommand{\mcc}[1]{{}}

\numberwithin{equation}{section}

\title[Semisimplicity of Lyapunov spectrum of irreducible cocycles]
{Semisimplicity of the Lyapunov spectrum for irreducible cocycles}

\author{Alex Eskin}
\thanks{Research  of  the first author is partially supported  by
NSF grants DMS 0244542, DMS 0604251 and DMS 0905912}
\address{
Department of Mathematics,
University of Chicago,
Chicago, Illinois 60637, USA\\
}
\email{eskin@math.uchicago.edu}

\author{Carlos Matheus}
\thanks{Research of the second author is partially supported by the Balzan project of Jacob Palis and by the French ANR grant
``GeoDyM'' (ANR-11-BS01-0004).}
\address{CMLS, \'Ecole Polytechnique, CNRS (UMR 7640), 91128 Palaiseau, France \\}
\email{matheus.cmss@gmail.com}

\date{\today}

\begin{document}

\begin{abstract}
  Let $G$ be a semisimple Lie group acting on a space $X$, let $\mu$ be a
  symmetric compactly supported measure on $G$, and let $A$ be a strongly
  irreducible linear cocycle over the action of $G$. We then have a
  random walk on $X$, and let $T$ be the associated shift map. We show
  that, under certain assumptions, the cocycle $A$ over the action of $T$ is conjugate to a block
  conformal cocycle.

  This statement is used in the recent paper by Eskin-Mirzakhani on
  the classification of invariant measures for the $SL(2,\mathbb{R})$ action on
  moduli space. The ingredients of the proof are essentially contained
  in the papers of Guivarch and Raugi and also Goldsheid and Margulis.
\end{abstract}

\maketitle

\vspace{-3.1pt}
\setcounter{tocdepth}{1}
\tableofcontents


\section{Introduction}\label{s.introduction}

\subsection{Statement of the main results}

Let $G$ be a semisimple Lie group. Denote by $\mu$ a symmetric compactly supported probability measure on $G$.

Let $X$ be a space where $G$ acts and denote by $\nu$ a $\mu$-stationary measure (that is, $\mu\ast\nu=\nu$ where $\mu\ast\nu:=\int_{G} g_*\nu d\mu(g)$). We assume that $\nu$ is $\mu$-ergodic.

Consider $L$ a real finite-dimensional vector space and $A:G\times
X\to SL(L)$ a (linear) cocycle\footnote{I.e., the cocycle relation $A(g_2 g_1, x) = A(g_2, g_1(x)) \cdot A(g_1, x)$ holds for all $x\in X$ and $g_1, g_2\in G$.}. Since it is sufficient for our purposes, we will assume that $A(g,x)$ is bounded for $g$ in the
  support of $\mu$.
Denote by $\bf{H}$ the \emph{algebraic hull} of $A(.,.)$ in Zimmer's sense, that is, the smallest linear $\mathbb{R}$-algebraic subgroup\footnote{Recall that the algebraic hull is unique up to conjugation (cf. Zimmer's book \cite{Zimmer}).} $\bf{H}$ such that there exists a measurable map $C:X\to SL(L)$ with $C(g(x))A(g,x)C(x)^{-1}\in \bf{H}$ for $\mu$-almost all $g\in G$ and $\nu$-almost all $x\in X$. In what follows, we will assume that $\bf{H}$ is a $\mathbb{R}$-simple Lie group with finite center, and a basis of $L$ is (measurably) chosen at each $x\in X$ so that the cocycle $A(.,.)$ takes its values in $\bf{H}$.

\begin{definition} We say that the cocycle $A(.,.)$ has an \emph{invariant system of subspaces} if there are measurable families $W_i(x)$, $i=1,\dots, k$, of subspaces of $L$ such that $A(g,x)(W_1(x)\cup\dots \cup W_k(x))\subset W_1(g(x))\cup\dots \cup W_k(g(x))$ for $\mu$-almost every $g\in G$ and $\nu$-almost every $x\in X$.
\end{definition}

\begin{definition}[Strong irreducibility] We say that $A(.,.)$ is \emph{strongly irreducible} if there are no  non-trivial and proper invariant systems of subspaces.
\end{definition}

We will be interested in the behavior of a strongly irreducible cocycle $A(.,.)$ on the Lyapunov subspaces obtained after multiplying the matrices $A(g,x)$ while following a random walk on $G$. For this reason, let us introduce the following objects.

Let $\Omega=G^{\mathbb{N}}$. Denote by $T:\Omega\times X\to \Omega\times X$ the natural forward shift map on $\Omega\times X$:
$$T(u,x)=(\sigma(u),u_1(x))$$
where $\sigma(u)=(u_2,\dots)$ for $u=(u_1, u_2,\dots)\in\Omega$. Denoting by $\beta=\mu^{\mathbb{N}}$ the probability measure on $\Omega$ naturally induced by $\mu$, it follows from the fact that $\nu$ is $\mu$-stationary that the probability measure $\beta\times\nu$ is $T$-invariant.

As we already mentioned above, from now on, we will assume that the stationary measure $\nu$ is $\mu$-ergodic\footnote{By definition, $\nu$ is $\mu$-ergodic if it is not a non-trivial convex combination of two distinct $\mu$-stationary measures. The fact that $\nu$ is $\mu$-ergodic is equivalent to the $T$-ergodicity of $\beta\times\nu$ is classical: see e.g. Benoist--Quint's book \cite{BQb}.}, that is, $\beta\times\nu$ is $T$-ergodic.

In this language, we can study the products of matrices of the cocycle $A(.,.)$ along random walks with the aid of the cocycle dynamics $F_A:\Omega\times X\times \textbf{H} \to\Omega\times X\times \textbf{H}$ naturally associated to $A(.,.)$:
$$F_A(u,x,h)=(T(u,x), A(u_1,x)h)$$
Actually, for our purposes, the ``fiber dynamics'' of $F_A$ will be more important than the base dynamics $T$. For this reason, given $u\in\Omega$ and $x\in X$, let us denote by $A^n(u,x)$ the matrix given by the formula:
\begin{eqnarray*}
F_A^n(u,x,\textrm{Id})&=&(T^n(u,x), A(u_n, u_{n-1}\dots u_1(x))\dots A(u_1,x))\\ 
&=& (T^n(u,x), A(u_n\dots u_1, x)) = (T^n(u,x), A^n(u,x)).
\end{eqnarray*}

In this context, the multiplicative ergodic theorem of V. Oseledets \cite{Oseledec} says that, if $\int\log^+\|A(g,x)\|d\mu(g)d\nu(x)<\infty$, then
there is a collection of numbers $\lambda_1 > \dots> \lambda_k$
with multiplicities $m_1,\dots, m_k$ called \emph{Lyapunov exponents} and, at
$\beta\times\nu$-almost every point $(u,x)\in \Omega\times X$, we have a \emph{Lyapunov flag}
\begin{equation}
\label{e.forward-flag}
\{0\}=V_{k+1}^+ \subset V_k^+(u,x) \subset\dots\subset V_1^+(u,x)=L
\end{equation}
such that $V_i^+(u,x)$ has dimension $m_i+\dots+m_k
$ and
$\lim\limits_{n\to\infty}\frac{1}{n}\log\|A^n(u,x)\vec{p}\|=\lambda_i$
whenever $\vec{p}\in V_i^+(u,x)\setminus V_{i+1}^+(u,x)$.

In this paper, we will study the consequences of the strong irreducibility of a cocycle for its Lyapunov spectrum (i.e., collection of Lyapunov exponents and flags). In particular, we will focus on the following property:

\begin{definition} We say that $F_A$ or simply $A(.,.)$ has
  \emph{semisimple Lyapunov spectrum} if its algebraic hull
  $\textbf{H}$ is \emph{block-conformal} in the sense that, for each
  $i=1,\dots, k$, $V_i^+(u,x)/V_{i+1}^+(u,x)$ possesses
an invariant splitting,
\begin{displaymath}
V_i^+(u,x)/V_{i+1}^+(u,x) = \bigoplus_{j=1}^{n_i} E_{ij}(u,x),
\end{displaymath}
and on each $E_{ij}(u,x)$ there exists  a (non-degenerate)
  quadratic form $\langle.,.\rangle_{ij,u,x}$ such that, for all
  $\vec{p}, \vec{q}\in E_{ij}(u,x)$ and for all
  $n\in\mathbb{N}$,
$$\langle A^n(u,x)\vec{p},
A^n(u,x)\vec{q}\rangle_{ij,T^n(u,x)}=e^{\lambda_{ij}(u,x,n)}\langle\vec{p},\vec{q}\rangle_{ij,(u,x)}$$
for some cocycle\footnote{That is,
  $\lambda_{ij}(u,x,m+1)=\lambda_{ij}(T^{m-1}(u,x),1)+\lambda_{ij}(u,x,m)$.} $\lambda_{ij}:\Omega\times X\times\mathbb{N}\to\mathbb{R}$.
\end{definition}

\subsection*{Standing assumptions} From now on, besides the hypotheses
\begin{itemize}
\item[(A1)] $G$ is a semisimple Lie group group acting on a space $X$;
\item[(A2)] $\mu$ is a symmetric compactly supported probability measure on $G$ and $\nu$ is an ergodic $\mu$-stationary probability measure on $X$;
\item[(A3)] $A:G\times X\to SL(L)$ is a linear cocycle (where $L$ is a real finite-dimensional vector space) such that $A(g,x)$ is bounded for $g$ in the support of $\mu$;
\item[(A4)] the algebraic hull $\textbf{H}$ of $A(.,.)$ is a $\mathbb{R}$-simple Lie group with finite center.
\item[(A5)] $A$ verifies Oseledets' integrability condition $\int\log\|A(g,x)^{\pm1}\|d\mu(g)d\nu(x)<\infty$,
\end{itemize}
we will actually require the (stronger assumption of) invariance of $\nu$ under $\textrm{supp}(\mu)$:
\begin{itemize}
\item[(A6)] for all $g\in\textrm{supp}(\mu)$, one has $g_*\nu=\nu$. 
\end{itemize} 
In fact, most of the arguments in this paper need just\footnote{We emphasize this point by writing most of this article in the setting of stationary measures.} the $\mu$-stationarity of $\nu$: as it turns out, the invariance of $\nu$ under $\textrm{supp}(\mu)$ is used \emph{only} at Subsection \ref{sec:subsec:proof:osceledets:a}.

In this context, the main result of this paper is the following theorem:

\begin{theorem}
\label{theorem:semisimple:lyapunov}
If $A(.,.)$ is strongly irreducible, then it has semisimple Lyapunov spectrum.

Furthermore, the top Lyapunov exponent corresponds to a single
conformal block, that is, for $\beta\times\nu$-a.e. $(u,x)$ there are
a (non-degenerate) quadratic form $\langle.,.\rangle_{u,x}$
and a cocycle $\lambda:\Omega\times X\times\mathbb{N}\to\mathbb{R}$ such that
\begin{equation}
\label{e.top-single-block-forward}
\langle A^n(u,x)\vec{p}, A^n(u,x)\vec{q}\rangle_{T^n(u,x)}=e^{\lambda(u,x,n)}\langle\vec{p},\vec{q}\rangle_{u,x}
\end{equation}
for all $\vec{p}, \vec{q}\in V_1^+(u,x)/V_2^+(u,x)$.
\end{theorem}

In fact, the ingredients of the proof of this result are essentially contained in the articles of Goldsheid-Margulis \cite{Goldsheid:Margulis} and Guivarc'h-Raugi \cite{Guivarch:Raugi:Frontiere}, \cite{Guivarch:Raugi:Contraction}. In particular, the fact that such a result holds is no surprise to the experts.

Nevertheless, we decided to write down a proof of this theorem here mainly for two reasons: firstly, this precise statement is hard to locate in these references, and, secondly, this result is relevant in the recent paper \cite{Eskin:Mirzakhani:measures} where a Ratner-type theorem is shown for the action of $SL(2,\mathbb{R})$ on moduli spaces of Abelian differentials.

\subsection{The backwards cocycle}

As it turns out, for the application in Eskin-Mirzakhani paper \cite{Eskin:Mirzakhani:measures}, one needs the analog of Theorem \ref{theorem:semisimple:lyapunov} for the backward shift.

More precisely, let $\Omega^-=G^{\mathbb{Z}-\mathbb{N}}$ and $\hat{\Omega}=\Omega^-\times\Omega$. Denote by
$T^-:\Omega^-\times X\to\Omega^-\times X$ the natural backward shift map on $\Omega^-\times X$:
$$T^-(v,y)=(\sigma^-(v), v_0^{-1}(y))$$
where $\sigma^-(v)=(\dots, v_{-1})$ for $v=(\dots, v_0)\in\Sigma^-$. Similarly, denote by $\hat{T}:\hat{\Omega}\times X\to\hat{\Omega}\times X$ the natural forward shift map on $\hat{\Omega}\times X$:
$$\hat{T}(v,u,x)=(\hat{\sigma}(v,u),u_1(x))$$
where $\hat{\sigma}(v,u)=(c_{i-1})_{i\in\mathbb{Z}}$ for $(v,u)=(c_{i})_{i\in\mathbb{Z}}$.

Recall that $\Omega$ is equipped with the probability measure $\beta=\mu^{\mathbb{N}}$, so that $\beta\times\nu$ is a $T$-invariant probability measure on $\Omega\times X$. Note that Borel measures on $\Omega$ and $\hat{\Omega}$ are uniquely determined by their values on cylinders. In particular, the natural projection $\pi_+:\hat{\Omega}\times X\to\Omega\times X$ induces a bijection $(\pi_+)_*$ between the spaces of $\hat{T}$-invariant and $T$-invariant Borel probability measures, and, \emph{a fortiori}, there exists an unique probability measure $\widehat{\beta\times \nu}$ on $\hat{\Omega}\times X$ projecting to $\beta\times\nu$ under $(\pi_+)_*$. In this context, the natural $T^-$-invariant probability measure
$\beta^X$ constructed in Lemma 3.1 of Benoist and Quint \cite{Benoist:Quint} is $\beta^X=(\pi_-)_*\circ(\pi_+)_*^{-1}(\beta\times\nu):=(\widehat{\beta\times\nu})^-$, where $\pi_-:\Omega\times X\to\Omega^-\times X$ is the natural projection.

Similarly to the previous subsection, we can study the products of matrices of the cocycle $A(.,.)$ along backward random walks with the aid of the dynamical system $F_A^-:\Omega^-\times X\times \textbf{H} \to\Omega^-\times X\times \textbf{H}$ given by
$$F_A^-(v,y,h)=(T^-(v,y), A(v_0, v_0^{-1}(y))^{-1}h)$$
naturally associated to $A$, or, equivalently, the ``fiber'' dynamics $A^{-n}(v,y)$ given by the formula:
\begin{eqnarray*}
(F_A^-)^n(v,y,\textrm{Id})&=&((T^-)^n(v,y),A(v_{-(n-1)}, v_{-(n-1)}^{-1}\dots v_0^{-1}(y))^{-1}\dots A(v_0,v_0^{-1}(y))^{-1}) \\ &=:&((T^-)^n(v,y), A^{-n}(v,y))
\end{eqnarray*}

By Oseledets multiplicative ergodic theorem, if $\int\log^+\|A(g,x)^{\pm1}\|d\mu(g)d\nu(x)<\infty$, then we have a Lyapunov flag
\begin{equation}
\label{e.backward-flag}
\{0\}=V_{0}^-\subset V_1^-(v,y) \subset \dots\subset V^-_{k}(v,y)=L
\end{equation}
such that $V^-_j(v,y)$ has dimension $m_1+\dots+m_j$ and
$\lim\limits_{n\to\infty}\frac{1}{n}\log\|A^{-n}(v,y)\vec{q}\|=-\lambda_j$
for $\vec{q}\in V^-_j(v,y)\setminus V^-_{j-1}(v,y)$,
where $\lambda_i$ are the Lyapunov exponents of $F_A$ and $m_i$ are their multiplicities from the paragraph surrounding \eqref{e.forward-flag}.

In this setting, we will show the following:
\begin{theorem}
\label{theorem:semisimple:lyapunov:inverse} Suppose that $A(.,.)$ is strongly irreducible. Then, $F_A^-$ has semisimple Lyapunov spectrum.

Furthermore, the largest Lyapunov exponent corresponds to
a single conformal block, i.e., for
$\beta^X=(\widehat{\beta\times\nu})^-$-a.e. $(v,y)$ there are a
(non-degenerate) quadratic form
$\langle.,.\rangle_{v,y}$ and a
cocycle $\lambda:\Omega^-\times
X\times\mathbb{N}\to\mathbb{R}$
such that
\begin{equation}
\label{e.top-single-block-backward}
\langle A^{-n}(v,y)\vec{p},
A^{-n}(v,y)\vec{q}\rangle_{(T^-)^n(v,y)}=e^{\lambda(v,y,n)}\langle\vec{p},\vec{q}\rangle_{v,y}
\end{equation}
for all $\vec{p}, \vec{q}\in V_1^-(v,y)$.
\end{theorem}

\subsection{The invertible cocycle}

Both Theorem~\ref{theorem:semisimple:lyapunov} and
Theorem~\ref{theorem:semisimple:lyapunov:inverse} are derived as a
consequence of a theorem about the two-sided walk.
By Oseledets theorem applied to $\hat{T}$, the flags
\eqref{e.forward-flag} and \eqref{e.backward-flag} exist for
$\widehat{\beta\times\nu}$-a.e. $(v,u,x)\in\hat{\Omega}\times X$
(and, moreover, $V_i^+(v,u,x)=V_i^+(u,x)$ and
$V_j^-(v,u,x)=V_j^-(v,x)$).

Then for
$\widehat{\beta\times\nu}$-a.e. $(v,u,x)\in\hat{\Omega}\times X$,
let us define
$$\mathcal{V}_\ell(v,u,x)=V_{\ell}^+(u,x)\cap V_{\ell}^-(v,x)$$
for every $1\leq\ell\leq k$. By \cite[Lemma~1.5]{Goldsheid:Margulis},
$\mathcal{V}_{\ell}(v,u,x)$ has dimension $m_{\ell}$ and
$$V_i^+(u,x)=\bigoplus\limits_{\ell=i}^k \mathcal{V}_{\ell}(v,u,x)
\quad \textrm{ and } \quad V_j^-(u,x)=\bigoplus\limits_{\ell=1}^j
\mathcal{V}_{\ell}(v,u,x)$$

In particular, for $\widehat{\beta\times\nu}$-a.e. $(v,u,x)$,
$$V_j^+(u,x)/V_{j+1}^+(u,x) \simeq \mathcal{V}_j(v,u,x) \simeq
V_{j}^-(v,x)/V_{j-1}^-(v,x)$$

Using this information, we show below that Theorems \ref{theorem:semisimple:lyapunov} and  
\ref{theorem:semisimple:lyapunov:inverse} follow from the corresponding result for the two-sided walk:

\begin{theorem}
\label{theorem:semisimple:twosided:lyapunov}
If $A(.,.)$ is strongly irreducible, then it has semisimple Lyapunov
spectrum, in the sense that the restriction of $A^n(v,u,x)$ to each
$\cV_i(v,u,x)$ is block-conformal.

Furthermore, the top Lyapunov exponent corresponds to a single
conformal block, that is, for
$\widehat{\beta\times\nu}$-a.e. $(v,u,x)$ there are
a (non-degenerate) quadratic form
$\langle.,.\rangle_{v,u,x}$ and a
cocycle $\lambda:\hat{\Omega}\times
  X\times\mathbb{N}\to\mathbb{R}$ such that
\begin{equation}
\label{eq:top-single-block-twosided}
\langle A^n(v,u,x)\vec{p},
A^n(v,u,x)\vec{q}\rangle_{T^n(v,u,x)}=e^{\lambda(v,u,x,n)}\langle\vec{p},\vec{q}\rangle_{v,u,x}
\end{equation}
for all $\vec{p}, \vec{q}\in \cV_1(v,u,x)$.
\end{theorem}

\begin{remark} It is shown in
  \cite[Appendix~C]{Eskin:Mirzakhani:measures} that
if the algebraic hull $\textbf{H}$ is the
whole group $SL(L)$, then all Lyapunov exponents are associated to
single conformal blocks, i.e., for $\widehat{\beta\times\nu}$-a.e.
$(v,u,x)\in\Omega\times X$ and for each $1\leq i\leq k$, there are a
(non-degenerate) quadratic form $\langle.,.\rangle_{i,v,u,x}$ and a
cocycle $\lambda_i:\hat{\Omega}\times X\times\mathbb{N}\to\mathbb{R}$ such
that
$$\langle A^n(v,u,x)\vec{p}, A^n(v,u,x)\vec{q}\rangle_{i,\hat{T}^n(v,u,x)}=e^{\lambda_i(v,u,x,n)}\langle\vec{p},\vec{q}\rangle_{i,v,u,x}$$
for all $\vec{p}, \vec{q}\in \cV_i(v,u,x)$. Furthermore, analogous statements hold for the forward and backward
  walks.
\end{remark}

\begin{proof}[Proof of Theorems \ref{theorem:semisimple:lyapunov} and  
\ref{theorem:semisimple:lyapunov:inverse} assuming Theorem \ref{theorem:semisimple:twosided:lyapunov}] 
Denote by $\langle.,.\rangle_{ij,v,u,x}$ the inner-products coming from the block-conformality property ensured by Theorem \ref{theorem:semisimple:twosided:lyapunov}. We will show that for $(\beta\times\nu)$-almost every $(u,x)$, resp. $\beta^X$-almost every $(v,x)$, the conformal class of $\langle.,.\rangle_{ij,v,u,x}$ does not depend on $v$, resp. $u$ (and this will suffice to obtain Theorems \ref{theorem:semisimple:lyapunov} and \ref{theorem:semisimple:lyapunov:inverse}).

Given $\varepsilon>0$, we can select a compact subset $K\subset \hat{\Omega}\times X$ with $\widehat{\beta\times\nu}(K)>1-\varepsilon$ such that the functions $(v,u,x)\mapsto \langle.,.\rangle_{ij,v,u,x}$ are uniformly continuous on $K$. By ergodicity, if we take $0<\varepsilon<1/2$ and we consider the corresponding compact subset $K$ just described, it follows that there exists $Y\subset\hat{\Omega}\times X$ with $\widehat{\beta\times\nu}(Y)=1$ such that, for any $(v,u,x)\in Y$, the elements of the orbit $(\hat{T}^n(v,u,x))_{n\in\mathbb{Z}}$ belong to $K$ for a set of integers $n$ with asymptotic\footnote{By definition, $\mathcal{R}\subset\mathbb{Z}$ has asymptotic density $>\delta$ when  $\liminf\limits_{m\to+\infty}\frac{1}{m}\#\{n\in \mathcal{R}: 0 \leq n < m\}>\delta$ and $\liminf\limits_{m\to+\infty}\frac{1}{m}\#\{n\in \mathcal{R}: -m < n \leq 0\}>\delta$.} density $>1/2$. 

Next, we define 
$$[\vec{p},\vec{q}]_{ij,v,u,x}:=\frac{\langle \vec{p}, \vec{q} \rangle_{ij,v,u,x}}{\langle \vec{p}, \vec{p} \rangle_{ij,v,u,x}^{1/2} \langle \vec{q}, \vec{q} \rangle_{ij,v,u,x}^{1/2}}$$ 

Let $(v,u,x), (v',u,x)\in Y$, resp. $(v,u,x), (v,u',x)\in Y$. By Theorem \ref{theorem:semisimple:twosided:lyapunov}, 
\begin{eqnarray}\label{e.conformality}
& &[\vec{p}_n, \vec{q}_n]_{ij,\hat{T}^n(v,u,x)} = [\vec{p}, \vec{q}]_{ij,v,u,x}, \quad [\vec{p}_n, \vec{q}_n]_{ij,\hat{T}^n(v',u,x)} = [\vec{p}, \vec{q}]_{ij,v',u,x} \,\,\, \forall n\geq 0, \\ & & [\vec{p}_n, \vec{q}_n]_{ij,\hat{T}^n(v,u',x)} = [\vec{p}, \vec{q}]_{ij,v,u',x} \,\,\, \forall n\leq 0 \nonumber 
\end{eqnarray}
where $\vec{p}_n:=A^n(u,x)\vec{p}$ and $\vec{q}_n:= A^n(u,x)\vec{q}$ for all $n\geq 0$, resp. $\vec{p}_n:=A^{-n}(v,x)\vec{p}$ and $\vec{q}_n:= A^{-n}(v,x)\vec{q}$ for all $n\leq 0$.  

By construction of $Y$, we can select a subsequence $n_k\to+\infty$, resp. $n_k\to-\infty$ such that $\hat{T}^{n_k}(v,u,x), \hat{T}^{n_k}(v',u,x)\in K$, resp. $\hat{T}^{n_k}(v,u,x), \hat{T}^{n_k}(v',u,x)\in K$ $\forall k\in\mathbb{N}$. 

Since the points $\hat{T}^n(v,u,x)$ and $\hat{T}^n(v',u,x)$, resp. $\hat{T}^n(v,u,x)$ and $\hat{T}^n(v,u',x)$, approach each other as $n\to+\infty$, resp. $n\to-\infty$, it follows from the definition of $K$ (and our choice of $(n_k)_{k\in\mathbb{N}}$) that 
$$[\vec{p}_{n_k},\vec{q}_{n_k}]_{ij,\hat{T}^{n_k}(v,u,x)} - [\vec{p}_{n_k},\vec{q}_{n_k}]_{ij,\hat{T}^{n_k}(v',u,x)}\to 0,$$
resp. 
$$[\vec{p}_{n_k},\vec{q}_{n_k}]_{ij,\hat{T}^{n_k}(v,u,x)} - [\vec{p}_{n_k},\vec{q}_{n_k}]_{ij,\hat{T}^{n_k}(v,u,x)}\to 0$$
as $k\to\infty$. 

By plugging this into \eqref{e.conformality}, we see that 
$$[\vec{p}, \vec{q}]_{ij,v',u,x} = [\vec{p}, \vec{q}]_{ij,v,u,x}, \quad \textrm{resp. } [\vec{p}, \vec{q}]_{ij,v,u',x} = [\vec{p}, \vec{q}]_{ij,v,u,x},$$
whenever $(v,u,x), (v',u,x)\in Y$, resp. $(v,u,x), (v,u',x)\in Y$. 

In other terms, 
\begin{eqnarray}\label{e.conformality-sumanifolds} & & \langle\vec{p}, \vec{q}\rangle_{ij,v',u,x} = c(v',v,u,x)\langle\vec{p}, \vec{q}\rangle_{ij,v,u,x}, \quad\textrm{resp. } \\  & & \langle\vec{p}, \vec{q}\rangle_{ij,v,u',x} = c(v,u',u,x)\langle\vec{p}, \vec{q}\rangle_{ij,v,u,x}, \nonumber
\end{eqnarray}
whenever $(v,u,x), (v',u,x)\in Y$, resp. $(v,u,x), (v,u',x)\in Y$.  

On the other hand, for $\beta\times\nu$, resp. $\beta^X$ almost every $(u,x)$, resp. $(v,x)$, we can Borel measurably select $v=v(u,x)\in\Omega^-$, resp. $u=u(v,x)\in\Omega$ such that $(v,u,x)\in Y$ (because of von Neumann selection theorem, see Theorem A.9 at page 196 of Zimmer's book \cite{Zimmer}). By setting $\langle.,.\rangle_{u,x}:=\langle.,.\rangle_{v(u,x),u,x}, \textrm{ resp. } \langle.,.\rangle_{v,x}:=\langle.,.\rangle_{v, u(v,x),x},$ 
we obtain that the conclusions of Theorems \ref{theorem:semisimple:lyapunov} and  
\ref{theorem:semisimple:lyapunov:inverse} are valid for these choices of inner-products thanks to Theorem \ref{theorem:semisimple:twosided:lyapunov} and the conformality relations \eqref{e.conformality-sumanifolds}.
\end{proof}

The remainder of this paper is devoted to the proof of Theorem \ref{theorem:semisimple:twosided:lyapunov}.

\section{$\reals$-simple Lie groups}

Let $\bH$ be a $\reals$-simple Lie
group. We will always assume that $\bH$ is a linear algebraic group
with finite center.
Let $\theta$ denote a Cartan involution of $\bH$, and let $\bK$
denote the set of fixed points of $\theta$. Then,
$\bK$ is a maximal compact subgroup of $\bH$.

Let $\bA$ denote a
maximal $\reals$-split torus of $\bH$ such that $\theta(\bA) = \bA$,
and let $\Sigma$ denote the
associated root system. Let $\Sigma^+$ denote the set of positive
roots, and let $\Delta$ denote the set of simple roots.
Let $\bB$ denote the Borel subgroup of $\bH$ corresponding to
$\Sigma^+$. Let $W$ denote the Weyl group of $(\bH,\bA)$.

Let $\bA_+$ be the positive Weyl chamber, i.e.,
\begin{displaymath}
\bA_+ := \{ a \in \bA \st \alpha(\log a) \ge 0 \text{ for all } \alpha \in
\Sigma^+ \}.
\end{displaymath}
We have the
decomposition
\begin{equation}
\label{eq:KAK}
\bH = \bK \bA_+ \bK.
\end{equation}
If $g \in \bH$ is written as $g = k_1 a k_2$ where $k_1, k_2 \in \bK$
and $a \in \bA_+$, we write for $\alpha \in \Sigma^+$
\begin{equation}
\label{eq:def:alpha:g}
\alpha(g) = \alpha(\log a).
\end{equation}

We also have the Bruhat decomposition
\begin{displaymath}
\bH = \bigsqcup_{w \in W} \bB w \bB.
\end{displaymath}
Let $w_0 \in W$ be the longest root. Then, $\bB w_0 \bB$ is open and
dense in $\bH$. Let
\begin{equation}
\label{eq:def:J}
\text{$J \subset \bH/\bB$ denote the complement of $\bB
w_0 \bB/\bB$ in $\bH/\bB$. }
\end{equation}

Given a subset $I \subset \Delta$, let $\bP_I$ denote the parabolic
subgroup of $\bH$ associated\footnote{I.e., $\bP_I$ contains $\bA$ and its root system $(\bP_I, A)$ is $\Sigma^+\cup \Sigma_I$ where $\Sigma_I\subset \Sigma$ consists of roots whose expansions relative to $\Delta$ have vanishing coefficients at elements of $I$.} to $I$. We have the Langlands
decomposition
\begin{displaymath}
\bP_I = \bM_I \bA_I \bN_I,
\end{displaymath}
where
\begin{displaymath}
\bA_I = \{ a \in \bA \st \alpha(\log a) = 0  \text{ for all } \alpha \in I \}.
\end{displaymath}
The group $\bM_I$ is semisimple, and commutes with $\bA_I$. The group
$\bN_I$ is unipotent, and $\bN_I \lhd \bP_I$.

For later use, we denote $\bar{\bN}_I = w_0 \bN_I w_0^{-1}$ and let $J_I$ be the complement of $(\bB w_0 \bP_I)/\bP_I$ in $\bH/\bP_I$.

We will use the rest of this section to deduce some \emph{general} properties of the actions of elements of $\bH$ on $\bH/\bP_I$.  In particular, even though these properties help in the proof of Theorem \ref{theorem:semisimple:twosided:lyapunov}, we decided to present them in their own section because they have nothing to do with the cocycle $A$ but only with the group $\bH$.

\subsection{A lemma of Furstenberg}

\begin{definition}[{\bf $(\epsilon,\delta)$-regular}]
\label{def:epsilon:regular}
Suppose $\epsilon > 0$ and $\delta > 0$ are fixed.
A measure $\eta$ on $\bH/\bB$ is {\em
  $(\epsilon,\delta)$-regular} if for
any $g \in \bH$,
\begin{displaymath}
\eta(\Nbhd_\epsilon(gJ)) < \delta,
\end{displaymath}
where $J$ is as in (\ref{eq:def:J}).
A measure $\eta_I$ on $\bH/\bP_I$ is {\em
  $(\epsilon,\delta)$-regular} if for any $g \in \bH$,
\begin{displaymath}
\eta_I(\Nbhd_\epsilon(gJ_I)) < \delta,
\end{displaymath}
where $J_I$ is the complement of $(\bB w_0 \bP_I)/\bP_I$ in
$\bH/\bP_I$.
\end{definition}

\begin{lemma}[Furstenberg\footnote{Compare with \cite[Theorems 8.3 and 8.4]{Furstenberg}.}]
\label{lemma:strongly:regular}
Suppose $I \subset \Delta$,
$g_n \in \bH$ is a sequence
, and
$\eta_n$ is a sequence of uniformly $(\epsilon,\delta)$-regular
measures on $\bH/\bP_I$. Suppose $\delta
  \ll 1$.  Write
\begin{displaymath}
g_n = k_n a_n k_n',
\end{displaymath}
where $k_n \in \bK$, $k_n' \in \bK$ and $a_n \in \bA_+$.
\begin{itemize}
\item[{\rm (a)}] Suppose $I \subset \Delta$ is such that
for all $\alpha \in \Delta\setminus I$,
\begin{equation}
\label{eq:dk:dkplus1}
\alpha(a_n)  \to \infty.
\end{equation}
Then, for any subsequential limit $\lambda$ of $g_n \eta_n$, we have
\begin{equation}
\label{eq:tmp:K(n):Span:e1es}
k_n \bP_I  \to k_{\infty}\, \bP_I \quad and \quad \lambda(\{ k_{\infty} \, \bP_I \}) \ge 1-\delta
\end{equation}
for some element $k_{\infty}\in\bK$.
\item[{\rm (b)}] Suppose  $g_n \eta_n \to \lambda$ where $\lambda$ is
  some measure on $\bH/\bP_I$. Suppose also that there exists
an element $k_{\infty}$ such that $\lambda(\{k_{\infty} \bP_I \}) > 5 \delta$.
Then, as $n \to \infty$,
(\ref{eq:dk:dkplus1}) holds for all $\alpha\in\Delta \setminus I$. As a consequence, by part (a),
(\ref{eq:tmp:K(n):Span:e1es}) holds and  $\lambda(\{k_{\infty} \bP_I \}) \ge
1-\delta$. \mcc{check this proof}
\end{itemize}
\end{lemma}

\bold{Proof of (a).} Without loss of
generality, $k_n'$ is the identity (or else we replace $\eta_n$ by
$k_n' \eta_n$).

Let $\bar{\bN}_I = w_0 \bN_I w_0^{-1}$.
By our assumption \eqref{eq:dk:dkplus1}, for $\bar{n} \in \bar{\bN}_I$,
\begin{displaymath}
a_n \bar{n} \bP_I = (a_n \bar{n} a_n^{-1}) \bP_I \to \bP_I  \quad\text{ in
  $\bH/\bP_I$. }
\end{displaymath}
For any $z \in \bH/\bP_I$ such that $z \not\in J_I$, we may
write $z = \bar{n} \bP_I$ for some $\bar{n} \in \bar{\bN}_I$. Therefore,
$d(g_n z, k_n \bP_I) \to 0$, where $d( \cdot, \cdot)$ denotes some
distance on $\bH/\bP_I$.
It then follows from the $(\epsilon,\delta)$-regularity of
$\eta_n$ that (\ref{eq:tmp:K(n):Span:e1es}) holds, and
any limit of $g_n \eta_n$ must give weight at least
$1-\delta$ to  $k_{\infty}\, \bP_I$ (where $k_{\infty}$ is a subsequential limit of $k_n$).

\bold{Proof of (b).} This is similar to
\cite[Lemma 3.9]{Goldsheid:Margulis}.  There is a subsequence of the $g_n$
(which we again denote by $g_n=k_n a_n k_n'$) such that for all $\gamma \in
\Delta$, either $\gamma(a_n) \to \infty$ or $\gamma(a_n)$ is bounded.
After passing again
to a subsequence, we may assume that $k_n \to k_{\infty}$. Also, without loss
of generality, we may assume that $k_n'$ is the identity (or else we
replace $\eta_n$ by $k_n' \eta_n$).

Suppose there exists $\alpha \in
\Delta \setminus I$ such that (\ref{eq:dk:dkplus1}) fails. Let $I' \subset \Delta$ denote the set of $\gamma \in \Delta$ such that,
for $\gamma \in \Delta\setminus I'$, $\gamma(a_n) \to \infty$. Since we
are assuming that $\alpha \in \Delta\setminus I$ and $\alpha \not\in
\Delta\setminus I'$, we have $\Delta\setminus I \not\subset
\Delta\setminus I'$, and thus $I' \not\subset I$.

Let $\bar{\bN}_\alpha \subset \bar{\bN}$ denote the subgroup obtained by
exponentiating the root subspace $-\alpha$. We may write $\bar{\bN}_I
= \bar{\bN}_\alpha \bar{\bN}'$ for some subgroup $\bar{\bN}'$ of
$\bar{\bN}$. Note that the action by left multiplication by $g_n$ on
$\bH/\bP_I$ does not shrink the direction $\bar{\bN}_\alpha$.

Write\footnote{The $(\epsilon,\delta)$-regularity of $\eta_n$ and our assumption $\lambda(\{k_{\infty} \bP_I \}) > 5 \delta$ imply that $k_{\infty} \bP_I\notin J_I$. Hence, it is possible to write $k_{\infty} \bP_I$ as claimed.} $k_{\infty} \bP_I= \bar{n}_\alpha \bar{n}'\bP_I$, where $\bar{n}_\alpha \in
\bar{\bN}_\alpha$, $\bar{n}' \in \bar{\bN}'$. Then, for $z \in \bH/\bP_I$,
$g_n z$ does not converge to $k_{\infty} \, \bP_I$ unless $z \in
\bar{n}_\alpha \bar{\bN}' \bP_I$ or $z \in J_I$. In particular, since
$\bar{n}_\alpha \bar{\bN}' \bP_I \subset \bar{n}_\alpha J_I$ (because
$w_0 \bar{\bN}' w_0^{-1} \in  \bB w_\alpha w_0 \bB$), we obtain that if $g_nz$ converges to $k_{\infty}\,\bP_I$ then $z\in J_I\cup \bar{n}_\alpha J_I$.

On the other hand, since $\eta_n$ is $(\epsilon,\delta)$-regular,
\begin{displaymath}
\eta_n(\Nbhd_\epsilon(J_I \cup \bar{n}_\alpha J_I)) < 2\delta.
\end{displaymath}

Therefore $\lambda(k_{\infty} \, \bP_I) < 3\delta$
which is a contradiction. Thus $\alpha(g_n) \to \infty$ for all
$\alpha \in \Delta\setminus I$.
Now, by part (a),
(\ref{eq:tmp:K(n):Span:e1es}) holds, and $\lambda(k_{\infty} \,\bP_I) \ge 1-\delta$.
\qed\medskip

\subsection{The functions $\xi_\alpha(\cdot,\cdot)$ and
  $\hat{\sigma}_\alpha(\cdot, \cdot)$.}\label{subsec:alpha:sigma:hat}
Let $\omega_\alpha$ be the fundamental weight corresponding to
$\alpha$, i.e. for $\gamma \in \Delta$,
\begin{displaymath}
\langle \omega_\alpha, \gamma \rangle = \begin{cases} 1 & \text{if
    $\alpha = \gamma$} \\
 0 & \text{if $\gamma \in \Delta\setminus \{\alpha\}$.}
\end{cases}
\end{displaymath}
Then,
\begin{equation}
\label{eq:roots:weights}
\alpha = \sum_{\gamma \in \Delta} \langle \alpha, \gamma\rangle \omega_\gamma.
\end{equation}
We write
\begin{equation}
\label{eq:def:omega:alpha}
\omega_\alpha(g) = \omega_\alpha(\log a), \qquad\text{where $g = k_1 a k_2$, $k_1, k_2 \in \bK$, $a \in \bA_+$.}
\end{equation}
Note that for all $\alpha \in \Delta$ and all $g \in \bH$,
\begin{equation}
\label{eq:alpha:in:terms:of:omega:alpha}
\alpha(g) = \sum_{\gamma \in \Delta} \langle \alpha, \gamma\rangle
\omega_\gamma(g).
\end{equation}
\begin{lemma}
\label{lemma:properties:w:alpha}
For all $g_1 \in \bH$, $g_2\in \bH$, and for all $\alpha \in \Delta$,
\begin{equation}
\label{eq:omega:alpha:subadditive}
\omega_\alpha(g_1 g_2) \le \omega_\alpha(g_1) + \omega_\alpha(g_2).
\end{equation}
and
\begin{equation}
\label{eq:omega:alpha:lower:subadditive}
\omega_\alpha(g_1 g_2) \ge \omega_\alpha(g_1) - \omega_\alpha(g_2^{-1}).
\end{equation}
\end{lemma}

\bold{Proof.} There exists a representation $\rho_\alpha: \bH \to
GL(V)$
such that
its highest weight is $\omega_\alpha$ (see \cite[Chapter V]{Knapp}). Let $\| \cdot \|$ be any
$\bK$-invariant norm on $V$. Then, since $\omega_\alpha$ is the
highest weight,
\begin{displaymath}
\| \rho_\alpha(g) \| \equiv \sup_{v \in V\setminus \{0\}} \frac{
  \|\rho_\alpha(g) v
  \|}{\|v\|} =
e^{\omega_\alpha(g)}.
\end{displaymath}
Since $\|\rho_\alpha(g_1 g_2) \| \le \|\rho_\alpha(g_1) \|
\|\rho_\alpha(g_2)\|$, (\ref{eq:omega:alpha:subadditive}) follows.

Now write $g_1 = h_1 h_2$, $g_2 = h_2^{-1}$, so that $g_1 g_2 =
h_1$. Substituting into (\ref{eq:omega:alpha:subadditive}), we get
\begin{displaymath}
\omega_\alpha(h_1) \le \omega_\alpha(h_1 h_2) + \omega_\alpha(h_2^{-1})
\end{displaymath}
which immediately implies (\ref{eq:omega:alpha:lower:subadditive}).
\qed\medskip

Let $\bP_\alpha$ be the parabolic subgroup corresponding to the subset
$\Delta \setminus \{\alpha \} \subset \Delta$. We can write
\begin{displaymath}
\bP_\alpha = \bM_\alpha \bA_\alpha \bN_\alpha,
\end{displaymath}
where
\begin{displaymath}
\bA_\alpha = \{ a \in \bA \st \gamma(\log a) = 0 \text{ for all } \gamma \in
\Delta\setminus \{\alpha\} \}.
\end{displaymath}
Note that $\bA_\alpha$ is one dimensional, and that $\bM_\alpha$
commutes with $\bA_\alpha$.
We have the Iwasawa decomposition
\begin{displaymath}
\bH =  \bK \bP_\alpha = \bK \bM_\alpha \bA_\alpha \bN_\alpha.
\end{displaymath}
Let $\bP_\alpha^{0} = \bM_\alpha \bN_\alpha$.
If we decompose $g \in \bH$ as $g = k_\alpha m_\alpha a_\alpha
n_\alpha$ with $k_\alpha \in \bK$, $m_\alpha \in
\bM_\alpha$, $a_\alpha \in \bA_\alpha$ and $n_\alpha \in \bN_\alpha$, then the
decomposition is unique up to the transformation $k_\alpha \to k_\alpha
m_1$, $m_\alpha \to
m_1^{-1} m_\alpha$ for $m_1 \in \bK \cap \bM_\alpha$.  We can thus define the
function $\xi_\alpha: \bH/\bP_\alpha^0 \to \reals$ by
\begin{displaymath}
\xi_\alpha(g) = \omega_\alpha(\log a),
\quad \text{ where $g = k m a n$, $k \in \bK$, $m \in
\bM_\alpha$, $a \in \bA_\alpha$ and $n \in \bN_\alpha$.}
\end{displaymath}
By definition, we have for $a \in \bA_\alpha$,
\begin{equation}
\label{eq:xi:alpha:a:equivariance}
\xi_\alpha(g a) = \xi_\alpha(g) + \xi_\alpha(a).
\end{equation}
We now define for $g \in \bH$, $z \in \bH/\bP_\alpha^0$,
\begin{displaymath}
\xi_\alpha(g,z) = \xi_\alpha(gz) - \xi_\alpha(z).
\end{displaymath}
Then, in view of (\ref{eq:xi:alpha:a:equivariance}), for $a \in
\bA_\alpha$, $\xi_\alpha(g,za) = \xi_\alpha(g,z)$. Thus, we may consider $\xi_\alpha( \cdot, \cdot)$
to be a function from $\bH \cross (\bH/\bP_\alpha)$ to
$\reals$.
\begin{lemma}
\label{lemma:xi:alpha:properties}
We have for all $\alpha \in \Delta$:
\begin{itemize}
\item[{\rm (a)}] For all $g_1, g_2 \in \bH$,
\begin{displaymath}
\xi_\alpha(g_1 g_2, z) = \xi_\alpha( g_1, g_2 z) + \xi_\alpha(g_2,z).
\end{displaymath}
\item[{\rm (b)}]
For all $g \in \bH$ and all $z \in \bH/\bP_\alpha$,
\begin{displaymath}
\xi_\alpha(g,z) \le \omega_\alpha(g),
\end{displaymath}
where  $\omega_\alpha(g)$ is as defined in (\ref{eq:def:omega:alpha}).
\item[{\rm (c)}] For all $\epsilon > 0$ there exists $C = C(\epsilon) > 0$
  such that for all $k_2 \in \bK$,
  for all $g \in \bK \bA_+ k_2$ and all $z \in \bH/\bP_\alpha$ with
  $d(k_2 z,
  J_{\alpha} ) > \epsilon$,
\begin{displaymath}
\xi_\alpha(g,z) \ge \omega_\alpha( g) - C.
\end{displaymath}
\end{itemize}
\end{lemma}

\bold{Proof.} Part (a) is clear from the definition of
$\xi_\alpha(\cdot, \cdot)$. To show part (b), note that there exists a
representation $\rho_\alpha: \bH \to GL(V)$ with highest weight
$\omega_\alpha$. Let $\| \cdot \|$ be any $\bK$-invariant norm on
$V$. Let $v_\alpha$ be the highest weight vector. Then $\bP_\alpha^0$ is the
stabilizer of $v_\alpha$, and for all $g \in \bH$,
\begin{displaymath}
\xi_\alpha(g) = \log \frac{\| \rho_\alpha(g) v_\alpha \|}{\|v_\alpha\|}.
\end{displaymath}
As in the proof of Lemma~\ref{lemma:properties:w:alpha},
\begin{displaymath}
\sup_{v \in V\setminus\{0\}} \log
\frac{\|\rho_\alpha(g) v\|}{\|v\|} = \omega_\alpha(g).
\end{displaymath}
Then, part (b) of Lemma~\ref{lemma:xi:alpha:properties} follows.

To show part (c), write $g = k_1 a k_2$, $k_1, k_2 \in \bK$, $a \in \bA_+$.
Note that if $d(k_2z, J_{\alpha})=d(k_2 z,
  (\bar{\bN}_\alpha \bP_\alpha)^c ) > \epsilon$,
then we can write
\begin{displaymath}
k_2 z = \bar{n}_\alpha \bP_\alpha,
\end{displaymath}
with $d(\bar{n}_\alpha,e) \le C_1(\epsilon)$.  Then,
$|\omega_\alpha(\bar{n}_\alpha)| < C(\epsilon)$. We have
\begin{align*}
\xi_\alpha(g,z) & = \xi_\alpha(k_1 a k_2, z) && \\
& = \xi_\alpha(k_1 a, k_2 z) && \text{by (a) and since
  $\xi_\alpha(k_2,z) = 0$ } \\
& = \xi_\alpha(k_1 a, \bar{n}_\alpha \bP_\alpha) && \\
& = \xi_\alpha(a, \bar{n}_\alpha \bP_\alpha) && \text{by (a) and since
  $\xi_\alpha(k_1, \cdot) = 0$ } \\
& = \xi_\alpha(a \bar{n}_\alpha, \bP_\alpha) -
\xi_\alpha(\bar{n}_\alpha,\bP_\alpha)  && \text{by (a)}\\
& \ge \xi_\alpha(a \bar{n}_\alpha, \bP_\alpha) - C(\epsilon)
&& \text{by (b) and since $|\omega_\alpha(\bar{n}_\alpha)| <
  C(\epsilon)$} \\
& = \xi_\alpha(a \bar{n}_\alpha a^{-1}, \bP_\alpha) + \xi_\alpha(a,
\bP_\alpha) - C(\epsilon) && \text{by (a)} \\
& =  \xi_\alpha(a \bar{n}_\alpha a^{-1}, \bP_\alpha) +
\omega_\alpha(a) - C(\epsilon) && \text{since
  $\xi_\alpha(a,\bP_\alpha) = \omega_\alpha(a)$} \\
& \ge \omega_\alpha(a) - 2 C(\epsilon) && \text{by (b) and since
  $|\omega_\alpha(a \bar{n}_\alpha a^{-1})| \le C(\epsilon)$.}
\end{align*}
\qed\medskip

For $\alpha \in \Delta$, $g \in \bH$, let $\bB=\bM\bA\bN$ be  Langlands' decomposition of $\bB$ and 
\begin{displaymath}
\hat{\sigma}_\alpha(g) = \alpha(a), \quad \text{ where $g = k m a
  n$, $k \in \bK$, $m \in
\bM$, $a \in \bA$ and $n \in \bN$.}
\end{displaymath}
Note that $\hat{\sigma}_\alpha$ descends to a well-defined function on $\bH/(\bM\bN)$. 

By definition, we have for $a \in \bA$,
\begin{equation}
\label{eq:sigma:alpha:a:equivariance}
\hat{\sigma}_\alpha(g a) = \hat{\sigma}_\alpha(g) + \hat{\sigma}_\alpha(a).
\end{equation}
We now define for $g \in \bH$, $z \in \bH/(\bM \bN)$,
\begin{displaymath}
\hat{\sigma}_\alpha(g,z) = \hat{\sigma}_\alpha(gz) - \hat{\sigma}_\alpha(z).
\end{displaymath}
Then, in view of (\ref{eq:sigma:alpha:a:equivariance}), for $a \in
\bA$, $\hat{\sigma}_\alpha(g,za) =
\hat{\sigma}_\alpha(g,z)$. Thus, we may consider $\hat{\sigma}_\alpha( \cdot, \cdot)$
to be a function $\bH \cross \bH/\bB \to \reals$.

\begin{lemma}
\label{lemma:sigma:alpha:properties}
We have for all $\alpha \in \Delta$:
\begin{itemize}
\item[{\rm (a)}] For all $g_1, g_2 \in \bH$ and $z \in \bH/\bB$,
\begin{displaymath}
\hat{\sigma}_\alpha(g_1 g_2, z) = \hat{\sigma}_\alpha( g_1, g_2 z) +
\hat{\sigma}_\alpha(g_2,z).
\end{displaymath}
\item[{\rm (b)}] For all $\epsilon > 0$ there exists $C = C(\epsilon) > 0$
  such that for all $k_2 \in \bK$,
  for all $g \in \bK \bA_+ k_2$ and all $z \in \bH/\bB$ with $d(k_2 z,
  J_{\alpha}) > \epsilon$,
\begin{displaymath}
\hat{\sigma}_\alpha(g,z) \ge \alpha( g) - C,
\end{displaymath}
where $\alpha(g)$ is as defined in (\ref{eq:def:alpha:g}).
\end{itemize}
\end{lemma}

\bold{Proof.} The natural map $\bH/\bP_\alpha \to
\bH/\bB$ allows us to consider the functions $\xi_\alpha(\cdot,
\cdot)$ to be functions $\bH \cross \bH/\bB \to \reals$. Then, in view
of (\ref{eq:roots:weights}), we have
\begin{displaymath}
\hat{\sigma}_\alpha(g,z) = \sum_{\gamma \in \Delta} \langle
\alpha,\gamma\rangle \xi_\gamma(g,z).
\end{displaymath}
Then (a) immediately follows from (a) of
Lemma~\ref{lemma:xi:alpha:properties}. Also,
\begin{align*}
\hat{\sigma}_\alpha(g,z) & = \langle \alpha,\alpha\rangle
\xi_\alpha(g,z) + \sum_{\gamma \ne \alpha} \langle \alpha,\gamma\rangle
\xi_\gamma(g,z) \\
& \ge \langle \alpha,\alpha\rangle
\xi_\alpha(g,z) + \sum_{\gamma \ne \alpha} \langle \alpha,\gamma\rangle
\omega_\gamma(g) && \text{by Lemma~\ref{lemma:xi:alpha:properties}(b)
  and since $\langle
  \alpha,\gamma\rangle \le 0$} \\
& \ge \langle \alpha,\alpha\rangle \omega_\alpha(g) - C(\epsilon) + \sum_{\gamma \ne \alpha} \langle \alpha,\gamma\rangle
\omega_\gamma(g) && \text{by Lemma~\ref{lemma:xi:alpha:properties}(c)}
\\
& = \alpha(g) - C(\epsilon).
\end{align*}
This completes the proof of (b).
\qed\medskip

\section{Cocycles with values in $\reals$-simple Lie groups}

Let $A:G\times X\to SL(L)$ be a linear cocycle satisfying the properties (A1) to (A6) described in \S\ref{s.introduction} above. In particular, we will assume that $A(.,.)$ takes values in its algebraic hull $\bH$. Furthermore, we will suppose that $\bH$ is a $\reals$-simple Lie group with finite center.

For $\alpha \in \Delta$, let
\begin{equation}
\label{eq:def:lambda:alpha}
\lambda_\alpha \equiv \limsup_{n \to +\infty} \frac{1}{n}
\alpha(A^n(u,x))
\end{equation}

By (\ref{eq:alpha:in:terms:of:omega:alpha}) and
Lemma~\ref{lemma:properties:w:alpha}, the map
\begin{displaymath}
(u,x,n) \to \lambda_\alpha(A^n(u,x))
\end{displaymath}
is a linear combination of subadditive cocycles.

Therefore, by the subadditive ergodic theorem, the limsup is
actually a limit. Also, by the ergodicity of $T$,
$\lambda_\alpha$ is constant a.e. on $\Omega \cross X$.

From now on, let us fix $I \subset \Delta$ minimal such that for $(\beta \cross \nu)$-a.e. 
$(u,x) \in \Omega \cross X$, we have
\begin{equation}
\label{eq:def:I}
I = \{ \alpha \in \Delta  \st  \lambda_\alpha = 0 \}.
\end{equation}
Thus, for all $\alpha \in \Delta\setminus I$, $\lambda_\alpha > 0$.

We will deduce Theorem~\ref{theorem:semisimple:twosided:lyapunov} from the
following:
\begin{theorem}
\label{theorem:goodroots}
Let $I \subset \Delta$ be as in (\ref{eq:def:I}). Then, for almost all
$((v,u),x) \in \hat{\Omega} \cross X$ there exists
$C(v,u,x) \in \bH$ such that 
\begin{equation}
\label{eq:conjugate:into:KM:AI}
C(\hat{T}^n(v,u,x))^{-1} A^n(v,u,x) C(v,u,x) = k_n(v,u,x) a_n(u,v,x),
\end{equation}
where $k_n(v,u,x) \in \bK \cap \bM_I$ and $a_n(v,u,x) \in \bA_I$,
and for all $\alpha \in \Delta\setminus I$,
\begin{equation}
\label{eq:theorem:goodroots}
\lim_{|n| \to \infty} \tfrac{1}{n} \alpha(\log
  a_n(v,u,x)) = \lambda_\alpha > 0
\end{equation}
where $\lambda_\alpha$ is as in (\ref{eq:def:lambda:alpha}).
\end{theorem}


Let $w_0 \in W$ be the longest root. Let $I' \subset \Delta$ be
defined by:
\begin{equation}
\label{eq:def:Iprime}
I' = \{ -w_0 \alpha w_0^{-1} \st \alpha \in I \}.
\end{equation}

Theorem~\ref{theorem:goodroots} will be deduced from the following
results.

\begin{proposition}
\label{prop:osceledets:GB}
Let $I' \subset \Delta$ be as in (\ref{eq:def:Iprime}). Then,
\begin{itemize}
\item[{\rm (a)}] For almost all $(u,x) \in \Omega \cross X$ there
  exists $C^+(u,x) \in \bH$ such that for all $n$ and almost all $(u,x)$,
\begin{displaymath}
C^+(T^n(u,x))^{-1} A^n(u,x) C^+(u,x) \in \bP_{I'}.
\end{displaymath}
\item[{\rm (b)}] For almost all $(v,x) \in \Omega^- \cross X$ there
  exists $C^-(v,x) \in \bH$ such that for all $n$ and almost all
  $(v,x)$,
\begin{displaymath}
C^-(T^{-n}(v,x))^{-1} A^{-n}(v,x) C^-(v,x) \in \bP_{I}.
\end{displaymath}
\item[{\rm (c)}] For almost all $(v,u,x) \in \hat{\Omega} \cross X$,
\begin{equation}
\label{eq:prop:osceledets:GB:transversal}
C^+(u,x)^{-1} C^-(v,x) \in \bP_{I'} \, w_0 \, \bP_{I}.
\end{equation}
\end{itemize}
\end{proposition}
We note that Proposition~\ref{prop:osceledets:GB} is similar in spirit to the geometrical versions of Osceledets multiplicative ergodic theorem in the literature (see the survey \cite{Filip}). The standard proofs (see e.g. \cite{Goldsheid:Margulis}) are
based on the subadditive ergodic theorem. We give a proof in
\S\ref{sec:proof:osceledets:GB} below based on the martingale
convergence theorem. Parts of this proof will be used again in the
proof of Theorem~\ref{theorem:goodroots} in
\S\ref{sec:proof:theorem:goodroots}.

\section{Proof of Proposition~\ref{prop:osceledets:GB}}
\label{sec:proof:osceledets:GB}

\subsection{A Zero One Law}
\label{sec:subsec:zero:one}



Let $\nu$ be an ergodic stationary measure on $X$.
Let $\hat{X} = X \cross \bH/\bB$. We then have an
action of $G$ on $\hat{X}$, by
\begin{displaymath}
g \cdot (x, z) = (gx, A(g,x) z).
\end{displaymath}
Let $\hat{\nu}$ be an ergodic $\mu$-stationary measure on $\hat{X}$
which
projects to $\nu$ under the natural map $\hat{X} \to X$.
Note there is always at least one such: one chooses $\hat{\nu}$ to be
an extreme point among the $\mu$-stationary measures which project to $\nu$. If
$\hat{\nu} = \hat{\nu}_1 + \hat{\nu}_2$ where the $\hat{\nu}_i$ are
$\mu$-stationary measures then $\nu = \pi_*(\hat{\nu}) =
\pi_*(\hat{\nu}_1) + \pi_*(\hat{\nu}_2)$. Since $\nu$ is
$\mu$-ergodic, this implies that $\pi_*(\hat{\nu}_1) = \nu$ or $\pi_*(\hat{\nu}_2) = \nu$,
hence the $\hat{\nu}_1$ or $\hat{\nu}_2$ also project to
$\nu$. Since $\hat{\nu}$ is an extreme point among such measures,
we must have $\hat{\nu}_1 = \nu$ or $\hat{\nu}_2 = \hat{\nu}$. This $\hat{\nu}$
is $\mu$-ergodic.

We may write
\begin{displaymath}
d\hat{\nu}(x,z) = d\nu(x) \, d\eta_x(z),
\end{displaymath}
where $\eta_x$ is a measure on $\bH/\bB$.

\begin{lemma}[cf. \protect{\cite[Lemma 4.2]{Goldsheid:Margulis}}, cf. \protect{\cite[Th\'eor\`eme~2.6]{Guivarch:Raugi:Frontiere}}, 
  cf. \protect{\cite[Lemma~C.10]{Eskin:Mirzakhani:measures}}]
\label{lemma:non:atomic:lyapunov}
For almost all $x \in X$ and any $g \in \bH$,
\begin{displaymath}
\eta_x( g J) = 0,
\end{displaymath}
where $J$ is defined in (\ref{eq:def:J}).
\end{lemma}

\bold{Proof.} Let $d$ be the smallest number such that
there exists
a subset $E \subset X$ with $\nu(E) > 0$
and for all $x \in E$ an irreducible algebraic subvariety \mcc{check
  goldsheid-margulis} $J_x \subset \bH/\bB$ of dimension $d$
with $\eta_x(J_x) > 0$.
For $x \in E$, let $\cS(x)$ denote the set of irreducible algebraic
subvarieties of $\bH/\bB$ of dimension $d$ such that for $Q \in
\cS(x)$, $\eta_x(Q) > 0$.

Note that for a.e. $x \in X$, for any $Q_1 \in \cS(x)$, $Q_2 \in \cS(x)$
with $Q_2 \ne Q_1$,
\begin{displaymath}
\eta_x(Q_1 \cap Q_2) = 0.
\end{displaymath}
(since $Q_1 \cap Q_2$ is an algebraic subvariety of dimension lower
than $d$). Thus
\begin{displaymath}
\sum_{Q \in \cS(x)} \eta_{x}(Q) \le 1.
\end{displaymath}
Therefore $\cS(x)$ is at most countable. Moreover, by setting
\begin{equation}
\label{eq:fmax}
f(x) = \max_{Q \in \cS(x)} \eta_x(Q)
\end{equation}
and $\mathcal{S}_{\max}(x):=\{Q\in \cS(x): \eta_x(Q)=f(x)\}$, we see that $\mathcal{S}_{\max}(x)$ is finite.

Consider the measurable subset $\mathcal{S}_{\max}=\{(x,z): x\in E, z\in \mathcal{S}_{\max}(x)\}\subset \bigcup\limits_{x\in E} (\{x\}\times \mathcal{S}(x))$. By definition, for each $x\in X$, the fiber
$\{z\in \bH/\bB : (x,z)\in \mathcal{S}_{\max}\}:=\cS_{\max}(x)$ of $\mathcal{S}_{\max}$ at $x$ is a finite set, and, in particular, $\mathcal{S}_{\max}(x)$ is a countable union of compact sets. By a result of Kallman (see, e.g., the statement of Theorem A.5 in Appendix A of Zimmer's book \cite{Zimmer}), we can find a Borel measurable section for the restriction to $\mathcal{S}_{\max}$ of the natural projection $\pi:\hat{X}\to X$. In other words, one has a Borel measurable map $X\ni x\mapsto Q^{(1)}_x\in\mathcal{S}_{\max}(x)$ whose graph $E_1:=\{(x,Q^{(1)}_x)\in\hat{X}: x\in X\}$ is a measurable subset of $\mathcal{S}_{\max}$. If $E_1=\mathcal{S}_{\max}$, we get that $\mathcal{S}_{\max}$ is the graph of a section of $\pi$. Otherwise, we apply once more Kallman's theorem to $\mathcal{S}_{\max}-E_1$ in order to obtain a measurable subset $E_2$ of $\mathcal{S}_{\max}-E_1$ given by the graph of a Borel measurable map $\pi(E_2):=\{y\in X: \#\mathcal{S}_{\max}(y)\geq 2\}\ni x\mapsto Q^{(2)}(x)\in\mathcal{S}_{\max}(x)$ that we extend (in a measurable way) to $X$ by setting $Q^{(2)}(x):=Q^{(1)}(x)$ whenever $\#\mathcal{S}_{\max}(x)=1$. Since the fibers $\mathcal{S}_{max}(x)$ of $\mathcal{S}_{\max}$ are finite sets, by iterating this procedure at most countably many times, we obtain a non-empty subset $Z\subset \mathbb{N}$ and, for each $m\in Z$, a Borel measurable map
$$X\ni x\mapsto Q^{(m)}(x)\in\mathcal{S}_{\max}(x)$$
such that $\mathcal{S}_{\max}(x)=\{Q^{(1)}(x),\dots, Q^{(\#\mathcal{S}_{\max}(x))}(x)\}$ for almost every $x\in X$.

Fix $m\in Z$. Since $\hat{\nu}$ is $\mu$-stationary, we have $\mu\ast\hat{\nu}(\textrm{graph}(Q^{(m)}))=\hat{\nu}(\textrm{graph}(Q^{(m)}))$, that is,

\begin{eqnarray}\label{eq:nu:hat:stat}
\int_X\eta_x(Q^{(m)}(x))d\nu(x) &= &\hat{\nu}(\textrm{graph}(Q^{(m)}))=\mu\ast\hat{\nu}(\textrm{graph}(Q^{(m)})) \\ &=&\int_G\int_{\hat{X}} \chi_{\textrm{graph}(Q^{(m)})}(gx,A(g,x)z)d\eta_x(z)d\nu(x)d\mu(g) \nonumber\\
& = & \int_G \int_X \eta_x(A(g,x)^{-1}Q^{(m)}(gx)) d\nu(x)d\mu(g) \nonumber\\
&=&  \int_X \left(\int_G \eta_x(A(g,x)^{-1}Q^{(m)}(gx)) d\mu(g)\right)d\nu(x). \nonumber
\end{eqnarray}

On the other hand, since $\eta_x(Q^{(m)}(x))=f(x)=\max_{Q \in \cS(x)} \eta_x(Q)$, we see that
\begin{equation}\label{eq:weight:trivial:estimate}
\int_G  \eta_x(A(g,x)^{-1}Q^{(m)}(gx)) d\mu(g)\leq f(x)=\eta_x(Q^{(m)}(x))
\end{equation}

By combining \eqref{eq:nu:hat:stat} and \eqref{eq:weight:trivial:estimate}, we deduce that
$$f(x)=\eta_x(Q^{(m)}(x))=\eta_x(A(g,x)^{-1}Q^{(m)}(gx)),$$
i.e., $A(g,x)^{-1}Q^{(m)}(gx)\in\mathcal{S}_{\max}(x)$ for $\mu$-almost every $g$ and $\nu$-almost every $x$. In other terms, for all $m\in Z$, $\mu$-almost every $g$ and $\nu$-almost every $x$, one has $Q^{(m)}(gx)\in A(g,x)\mathcal{S}_{\max}(x)$. By putting together this inclusion with the facts that $\mathcal{S}_{\max}(y)=\{Q^{(1)}(y),\dots, Q^{(\#\mathcal{S}_{\max}(y))}(y)\}$ for $\nu$-almost every $y$ and $\nu$ is $\mu$-stationary, one has that $\mathcal{S}_{\max}(gx)\subset A(g,x)\mathcal{S}_{\max}(x)$ for $\mu$-almost every $g$ and $\nu$-almost $x$.

Now, let $n_0\in Z$ the smallest integer in $Z$ such that $\{x\in X: \#\mathcal{S}_{\max}(x)\leq n_0\}$ has positive $\nu$-measure. Because $\mathcal{S}_{\max}(gx)\subset A(g,x)\mathcal{S}_{\max}(x)$ (for $\mu\times\nu$-almost every $(g,x)$), the set $\{x\in X: \mathcal{S}_{\max}(x)\leq n_0\}$ is essentially invariant. Thus, from the $\mu$-ergodicity of $\nu$ and our choice of $n_0$, we conclude that $\{x\in X: \#\mathcal{S}_{\max}(x)=n_0\}$ has full $\nu$-measure. Hence, from the $\mu$-stationarity of $\nu$, we obtain that $\#\mathcal{S}_{\max}(gx)=\#\mathcal{S}_{\max}(x)=n_0$ for $(\mu\times\nu)$-almost every $(g,x)$. In particular, the inclusion $\mathcal{S}_{\max}(gx)\subset A(g,x)\mathcal{S}_{\max}(x)$ is actually an equality 
\begin{equation}\label{eq:S-equiv}
\mathcal{S}_{\max}(gx)=A(g,x)\mathcal{S}_{\max}(x)
\end{equation} 
for $(\mu\times\nu)$-almost every $(g,x)$, that is, the cocycle $A$ permutes the finite sets $\cS_{\max}(x)$. 

Denote by $\cM$ the space of algebraic subvarieties of $\bH/\bB$ of dimension $d$, so that, by definition, $\cS_{\max}(x)$ are finite subsets of $\cM$. Note that the algebraic group $\bH$ acts  algebraically on $\cM$. By the Borel--Serre theorem (see, e.g., Theorem 3.1.3 of Zimmer's book \cite{Zimmer}), the action of $\bH$ on $\cM$ has locally closed orbits, and, hence, $\cM/\bH$ is Hausdorff. Consider the function $f:\Omega\times X\to\cM/\bH$, $f(u,x):=\bH\cS_{\max}(x)$. By \eqref{eq:S-equiv}, the function $f(u,x)$ is $T$-invariant. By the $T$-ergodicity of $\beta\times\nu$, it follows that $f$ is almost everywhere constant, so that there is $\cS_0\in\cM$ such that, for $\nu$-a.e. $x$, one has  $\cS_{\max}(x) = h(x) \cS_0$ with $h(x)\in \bH$. 

In particular, it follows that the conjugated cocycle $A'(g,x) = h(gx)^{-1} A(g,x) h(x)$ stabilizes $\cS_0$, i.e., $A'(g,x)$ takes its values in the stabilizer $\textrm{Stab}_{\bH}(\cS_0)$ of $\cS_0$ in $\bH$. Hence, the algebraic hull of $A(.,.)$ must be a subgroup of $\textrm{Stab}_{\bH}(\cS_0)$. This is a contradiction with our assumption that $\bH$ is the algebraic hull of $A(.,.)$ because $\textrm{Stab}_{\bH}(\cS_0)$ is a proper subgroup of $\bH$ (since $\bH$ acts transitively on $\bH/\bB$).
\qed\medskip

\subsection{Another lemma of Furstenberg}

Let $\hat{X} = X \cross \bH/\bB$.
The group $G$ acts on the space $\hat{X}$ is by
\begin{equation}
\label{eq:straight:action}
g \cdot (x,z) = (g x, A(g,x) z).
\end{equation}

We choose some ergodic $\mu$-stationary measure $\hat{\nu}$ on
$\hat{X}$, which projects to $\nu$, and write
\begin{displaymath}
d\hat{\nu}(x,z) = d\nu(x) \, d\eta_x(z).
\end{displaymath}
Note that Lemma~\ref{lemma:non:atomic:lyapunov} applies to the
measures $\eta_x$ on $\bH/\bB$.

\begin{lemma}[Furstenberg\footnote{Compare with \cite[Theorem 8.5]{Furstenberg}.}]
\label{lemma:furstenberg:formula}
For $\alpha \in \Delta$, let $\bar{\sigma}_\alpha:
G \cross \hat{X} \cross (\bH/\bB)\to \reals$
be given by
\begin{displaymath}
\bar{\sigma}_\alpha(g,x,z) = \hat{\sigma}_\alpha(A(g,x) z)
\end{displaymath}
with $\hat{\sigma}_{\alpha}(.,.)$ as in \S \ref{subsec:alpha:sigma:hat}.
Then,  we have
\begin{displaymath}
\lambda_\alpha = \int_{G}
\int_{\hat{X}} \bar{\sigma}_\alpha( g, x, z) \,
d\hat{\nu}(x,z) \, d\mu(g).
\end{displaymath}
where $\lambda_\alpha$ is as in (\ref{eq:def:lambda:alpha}).
\end{lemma}

\bold{Proof.} This is similar to the proof of \cite[Lemma
5.2]{Goldsheid:Margulis}. Note that
\begin{displaymath}
\xi_\alpha(g,z) = \sum_{\gamma \in \Delta} \langle \omega_\alpha,
\omega_\gamma\rangle \hat{\sigma}_\gamma(g,z)
\end{displaymath}
where $\xi_{\alpha}(.,.)$ and $\hat{\sigma}_{\alpha}(.,.)$ as in \S \ref{subsec:alpha:sigma:hat}.
Therefore, it is enough to show that for all $\alpha \in \Delta$,
\begin{equation}
\label{eq:enough:to:prove}
\int_{G} \int_{\hat{X}} \xi_\alpha( A(g, x) z) \,
d\hat{\nu}(x,z) \, d\mu(g) = \sum_{\gamma \in \Delta} \langle
\omega_\alpha,\omega_\gamma\rangle \lambda_\gamma.
\end{equation}
By (\ref{eq:def:lambda:alpha}), the boundedness assumption on
  $A(\cdot, \cdot)$ and the dominated convergence theorem, we have
\begin{displaymath}
\lambda_\alpha = \lim_{n \to \infty} \frac{1}{n} \int_{\Omega \cross X}
\alpha(A^n(g,x)) \, d\beta(g) d\nu(x).
\end{displaymath}
Thus,
\begin{displaymath}
\sum_{\gamma \in \Delta} \langle \omega_\alpha, \omega_\gamma\rangle
\lambda_\gamma = \lim_{n \to \infty} \frac{1}{n} \int_{\Omega \cross X}
\omega_\alpha(A^n(g,x)) \, d\beta(g) d\nu(x).
\end{displaymath}
Write $A^n(g,x) = \bar{k}_n(g,x) \bar{a}_n(g,x) \bar{k}'_n(g,x)$,
where $\bar{k}_n(g,x), \bar{k}_n'(g,x) \in \bK$, $\bar{a}_n(g,x) \in
\bA_+$, and fix $\epsilon > 0$.
Then, by Lemma~\ref{lemma:xi:alpha:properties}, for all
$z \in \bH/\bP_\alpha$ with
  $d(\bar{k}_n'(g,x) z,
  (\bar{\bN}_\alpha \bP_\alpha)^c ) > \epsilon$, we have
\begin{displaymath}
\omega_\alpha(A^n(g,x)) \ge \xi_\alpha(A^n(g,x),z) \ge
\omega_\alpha(A^n(g,x)) -C(\epsilon).
\end{displaymath}
Hence, by Lemma~\ref{lemma:non:atomic:lyapunov}, 
\begin{displaymath}
\sum_{\gamma \in \Delta} \langle \omega_\alpha, \omega_\gamma\rangle
\lambda_\gamma = \lim_{n \to \infty} \frac{1}{n} \int_{\Omega \cross \hat{X}}
\xi_\alpha(A^n(g,x), z) \, d\beta(g) \, d\hat{\nu}(x,z).
\end{displaymath}
By Lemma~\ref{lemma:xi:alpha:properties} (a) and by iterating the cocycle relation for $A(.,.)$,
\begin{displaymath}
\xi_\alpha(A^n(g,x),z) = \sum_{k=1}^n \xi_\alpha( A(g_k, g_{k-1} \dots
g_1 x), A(g_{k-1} \dots g_1,x) z),
\end{displaymath}
\mcc{check indices}
Since $\hat{\nu}$ is stationary, each of the terms in the sum has the
same integral over $\Omega \cross \hat{X}$ (with respect to $\beta \cross \hat{\nu}$). Therefore
\begin{displaymath}
\frac{1}{n} \int_{\Omega \cross \hat{X}}
\xi_\alpha(A^n(g,x), z) \, d\beta(g) \, d\hat{\nu}(x,z) = \int_{G} \int_{\hat{X}} \xi_\alpha( A(g, x), z) \,
d\hat{\nu}(x,z) \, d\mu(g),
\end{displaymath}
which completes the proof of (\ref{eq:enough:to:prove}).
\qed\medskip

\subsection{Proof of Proposition~\ref{prop:osceledets:GB}(a).}
\label{sec:subsec:proof:osceledets:a}
For $u \in \Omega$, let the measures $\nu_u$, $\hat{\nu}_u$ be essentially\footnote{It is shown in \cite[Lemma 3.2]{Benoist:Quint} that the convergence of $(u_1\dots u_n)_*\nu$ and $(u_1\dots u_n)_*\hat{\nu}$ for $\beta$-almost every $u=(u_1, u_2,\dots)\in\Omega$. In our setting, this implies the convergence of $(u_n\dots u_1)^{-1}_*\nu$ and $(u_n\dots u_1)^{-1}_*\hat{\nu}$ for $\beta$-almost $u\in\Omega$ because $(u_n\dots u_1)^{-1} = u_1^{-1}\dots u_n^{-1}$, $\beta=\mu^{\mathbb{N}}$, and $\mu$ is symmetric.} as
defined in \cite[Lemma 3.2]{Benoist:Quint}, i.e.
\begin{equation}\label{e.BQ1}
\nu_u = \lim_{n \to \infty} (u_n \dots u_1)^{-1}_* \nu
\end{equation}
\begin{equation}\label{e.BQ2}
\hat{\nu}_u = \lim_{n \to \infty} (u_n \dots u_1)^{-1}_*
\hat{\nu}.
\end{equation}
The limits exist for $\beta$-a.e. $u \in \Omega$ by the martingale convergence theorem.
We disintegrate
\begin{displaymath}
d\hat{\nu}(x,z) = d\nu(x) \, d\eta_x(z), \qquad
d\hat{\nu}_u(x,z) = d\nu_u(x) \, d\eta_{u,x}(z).
\end{displaymath}

We want to exploit  \eqref{eq:straight:action}, \eqref{e.BQ1} and \eqref{e.BQ2} to deduce that the sequence of conditional measures $A((u_n \dots u_1)^{-1}, u_n \dots u_1 x)\eta_{u_n \dots u_1 x}$ converges to the conditional measures $\eta_{u,x}$. For this sake, we shall use the invariance of $\nu$ under $\textrm{supp}(\mu)$. More concretely, let $C\subset X$ and $D\subset \bH/\bB$ be measurable subsets, and denote by $\chi_A$ the characteristic function of $A$. Since $\hat{\nu}$ is $\mu$-stationary, we have from  \eqref{eq:straight:action} that 
$$
\int_C \eta_x(D)d\nu(x) = \hat{\nu}(C\times D) = (\mu\ast\hat{\nu})(C\times D) = \int \chi_C(gy) A(g, y)_*\eta_y(D) d\nu(y) d\mu(g)
$$
By using the invariance of $\nu$ to rewrite the right-hand side of this equality, we get 
\begin{eqnarray*}
\int_C \eta_x(D)d\nu(x) &=& \int \chi_C(x) A(g, g^{-1}x)_*\eta_{g^{-1}x}(D) d\nu(x) d\mu(g) \\ 
&=& \int_C\left(\int_G A(g,g^{-1}x)_*\eta_{g^{-1}x}(D) d\mu(g)\right) d\nu(x)
\end{eqnarray*}
Because $C$ and $D$ are arbitrary, we deduce that 
$$\eta_x=\int_{G} A(g,g^{-1}x)_*\eta_{g^{-1}x} d\mu(g)$$
From this identity, the cocycle relation and the symmetry of $\mu$, we conclude that 
$$A((g_{n-1}\dots g_1)^{-1}, g_{n-1}\dots g_1 x)_*\eta_{g_{n-1}\dots g_1 x} = \int_G A((g_n\dots g_1)^{-1}, g_n\dots g_1 x)_*\eta_{g_n\dots g_1 x} d\mu(g_n),$$ 
so that $A((g_n\dots g_1)^{-1}, g_n\dots g_1 x)_*\eta_{g_n\dots g_1 x}$ is a martingale. Thus, the martingale convergence theorem and the uniqueness of Rokhlin disintegration imply that
\begin{displaymath}
\lim_{n \to \infty} A((u_n \dots u_1)^{-1}, u_n \dots u_1 x)
\eta_{u_n \dots u_1 x} = \eta_{u,x}.
\end{displaymath}
for $(u,x)$ in a set of $\beta \cross \nu$ full measure.

Note that, by the cocycle relation $A(g^{-1} ,g x) = A(g,x)^{-1}$, one has
\begin{displaymath}
A((u_n \dots u_1)^{-1}, u_n \dots u_1 x) = A(u_n \dots u_1,x)^{-1}.
\end{displaymath}
Hence, on a set of $\beta \cross \nu$-full measure,
\begin{equation}
\label{eq:measures:converge}
\lim_{n \to \infty} A(u_n \dots u_1,x)^{-1} \eta_{u_n \dots u_1 x} =
\eta_{u,x}.
\end{equation}

In view of Lemma~\ref{lemma:non:atomic:lyapunov} (see also the proof of
\cite[Lemma~14.4]{Eskin:Mirzakhani:measures}), given $\delta>0$, 
there exists a
compact  $\cK_{\delta} \subset X$ with $\nu(\cK_{\delta}) > 1-\delta$ and $\epsilon = \epsilon(\delta) > 0$ with $\epsilon(\delta)\to0$ as $\delta\to0$ such that the family of
measures $\{\eta_x \}_{x \in \cK_{\delta}}$ is uniformly
$(\epsilon,\delta/5)$-regular (in the sense of Definition \ref{def:epsilon:regular}).
Let
\begin{equation}\label{eq:Ndelta}
\cN_\delta(u,x) = \{ n \in \natls \st u_n \dots u_1 x \in
\cK_{\delta} \}.
\end{equation}
Write
\begin{equation}
\label{eq:decomp:Ainverse}
A( u_n \dots u_1, x)^{-1} = k_n(u,x) a_n(u,x) k_n'(u,x)
\end{equation}
where $k_n(u,x) \in \bK$, $k_n'(u,x) \in \bK$ and
$a_n(u,x) \in \bA_+$. We also
write
\begin{equation}
\label{eq:decomp:A:noinverse}
A( u_n \dots u_1, x) = \bar{k}_n(u,x) \bar{a}_n(u,x) \bar{k}_n'(u,x).
\end{equation}
where $\bar{k}_n$ and $\bar{k}_n'$ are elements of $\bK$, and
$\bar{a}_n \in \bA_+$.
Then, $\bar{a}_n(u,x) = w_0 a_n(u,x)^{-1} w_0^{-1}$ and thus,
\begin{multline}
\label{eq:d:dbar:conversion}
\alpha'(a_n(u,x)) = \alpha(\bar{a}_n(u,x)), \\
\bar{k}_n(u,x) = k'_n(u,x)^{-1} w_0^{-1}, \quad \bar{k}'_n(u,x) = w_0
k_n(u,x)^{-1},
\end{multline}
where $w_0$ is longest element of the Weyl group, and $\alpha' = -w_0
\alpha w_0^{-1}$.
 \mcc{check inverses}

Let $\pi_{I'}: \bH/\bB \to
\bH/\bP_{I'}$ be the natural map. Let $\eta_x^{I'} =
(\pi_{I'})_* \eta_x$ and
$\eta_{u,x}^{I'} = (\pi_{I'})_* \eta_{u,x}$. Then, $\eta_x^{I'}$ and
$\eta_{u,x}^{I'}$ are measures on $\bH/\bP_{I'}$.

Suppose $\alpha \in \Delta\setminus I$. Then, $\lambda_\alpha >
0$ and, \emph{a fortiori},
\begin{displaymath}
\lim_{n \to \infty} \alpha(\bar{a}_n(u,x)) \to \infty.
\end{displaymath}
Thus,
\begin{displaymath}
\lim_{n \to \infty} \alpha'(a_n(u,x)) \to \infty
\end{displaymath}
for each $\alpha'\in\Delta\setminus I'$.

Applying Lemma~\ref{lemma:strongly:regular}(a) to $g_n = A(u_n \dots
u_1,x)^{-1}$ for $n \in \cN_\delta(u,x)$ and the
$(\epsilon,\delta)$-regular measures $\eta_n =
\eta^{I'}_{u_n \dots u_1x}$ we get that there exists
$\bar{k} = \bar{k}(I',u,x) \in \bK$ such that, for $n \in \cN_\delta(u,x)$, one has $k_n(u,x) \bP_{I'} \to \bar{k} \bP_{I'}$ and
\begin{displaymath}
\eta^{I'}_{u,x}(\{ \bar{k}\, \bP_{I'}\} ) \ge 1- \delta.
\end{displaymath}
Since $\delta > 0$ is arbitrary, we get that for almost all $(u,x)$,
$\eta^{I'}_{u,x}$ is supported on one point of
$\bH/\bP_{I'}$. Now choose $C^+(u,x) \in \bH/\bB$ so that $\pi_{I'}(C^+(u,x)) = \bar{k}(I',u,x)
\bP_{I'}$. The desired property about $C^+(u,x)$ follows from the stationarity of $\hat{\nu}$.
\qed\medskip

\subsection{Proof of Proposition~\ref{prop:osceledets:GB} (b),(c).}
The proof of Proposition~\ref{prop:osceledets:GB} (b) is virtually
identical to the proof of Proposition~\ref{prop:osceledets:GB}(a), and
so we omit the details. Part (c) of
Proposition~\ref{prop:osceledets:GB} is also a classical fact,
cf. \cite[Lemma~1.5]{Goldsheid:Margulis}. We give an outline of a
geometric argument as follows.

Let $\bH/\bK$ be the symmetric space associated to $\bH$.
We say that two geodesic
rays (parametrized by arc length) are equivalent if they stay a
bounded distance apart.

By the geometric version of the multiplicative ergodic theorem
\cite{Karlsson:Margulis}, \cite{Kaimanovich},
for almost all $(u,x) \in
\Omega \cross X$
there exists a geodesic ray $\gamma^+: [0,\infty) \to \bH/\bK$
with $\gamma^+(0) = \bK$ such that
\begin{equation}
\label{eq:karlsson:margulis}
\lim_{n \to \infty} \frac{1}{n} d(A^n(u,x)^{-1} \bK, \gamma^+(n)) = 0.
\end{equation}
Similarly, by applying the same argument to the backwards walk, we get
that for almost all $(v,x) \in \Omega^- \cross X$
there exists a geodesic ray
$\gamma^-: [0,\infty) \to \bH/\bK$ such that
\begin{equation}
\label{eq:karlsson:margulis:backwards}
\lim_{n \to \infty} \frac{1}{n} d(A^{-n}(v,x)^{-1} \bK, \gamma^-(n)) =
0.
\end{equation}
Let $F = F(v,u,x)$ be a flat in $\bH/\bK$ which contains rays $\hat{\gamma}^+$
and $\hat{\gamma}^-$ equivalent to
$\gamma^+$ and $\gamma^-$ respectively. Then, we have
\begin{equation}
\label{eq:karlsson:margulis:hat}
\lim_{n \to \infty} \frac{1}{n} d(A^n(v,u,x)^{-1} \bK, \hat{\gamma}^+(n)) = 0.
\end{equation}
and
\begin{equation}
\label{eq:karlsson:margulis:backwards:hat}
\lim_{n \to \infty} \frac{1}{n} d(A^{-n}(v,u,x)^{-1} \bK, \hat{\gamma}^-(n)) = 0.
\end{equation}
Therefore, for every $\delta > 0$ there exists a set $K_\delta \subset
\hat{\Omega} \cross X$ with
$\widehat{\beta \cross \nu}(K_\delta) > 1-\delta$ and $N > 0$ such
that for $(v,u,x) \in K_\delta$ and $n > N$,
\begin{equation}
\label{eq:uniform:karlsson:margulis:hat}
d(A^n(v,u,x)^{-1} \bK, \hat{\gamma}^+(n)) \le \delta n, \text{ and }
d(A^{-n}(v,u,x)^{-1} \bK, \hat{\gamma}^-(n)) < \delta n.
\end{equation}
Let $X_n = A^n(v,u,x)^{-1} \bK$, and let $\hat{X}_n$ be the closest point
to $X_n$ on $\hat{\gamma}_n^+$. Then, by
(\ref{eq:karlsson:margulis:hat}), for $(v,u,x) \in K_\delta$ and
$n > N$,
\begin{equation}
\label{eq:d:Xn:hat:X:n}
d(X_n, \hat{X}_n) \le \delta n.
\end{equation}
Let $\hat{\gamma}^-_n(t)$ be unique
geodesic ray equivalent to $\hat{\gamma}^-$ such that $\hat{\gamma}^-_n(0) =
\hat{X}_n$. Then, as long as $T^n(v,u,x) \in K_\delta$, \mcc{explain}
and $m > N$, by (\ref{eq:uniform:karlsson:margulis:hat}) and
(\ref{eq:d:Xn:hat:X:n}), we have
\begin{displaymath}
d(A^{-m}(v,u,x)^{-1} X_n, \hat{\gamma}^-_n(m)) \le \delta n +
\delta m.
\end{displaymath}
Since $A^n(v,u,x)$ and $A^{-n}(\hat{T}^n(v,u,x))$ are inverses, we have
\begin{equation}
\label{eq:hat:gamma:minus:n:n:o:n}
d(\hat{\gamma}^-_n(n),e) \le 2 \delta n.
\end{equation}
Note that $\hat{X}_n$, $\hat{\gamma}^+$, $\hat{\gamma}_n^-$ all lie in
$F$. However in that case, (\ref{eq:hat:gamma:minus:n:n:o:n}) (for
sufficiently small $\delta$ and large enough $n$) implies
that
\begin{equation}
\label{eq:opposite:chambers}
\text{$\hat{\gamma}^+$ and $\hat{\gamma}^-$ belong to the closures of
  opposite Weyl chambers in $F$.}
\end{equation}
We now interpret (\ref{eq:opposite:chambers}) in terms of $C^+(u,x)$ and
$C^-(v,x)$. We can write
\begin{displaymath}
\gamma^+(t) = k(u,x) \hat{\Lambda}^t \bK,
\end{displaymath}
where $k(u,x) \in \bK$ and $\hat{\Lambda}^t \in \bA_+$. Then, by comparing
(\ref{eq:karlsson:margulis}) with (\ref{eq:decomp:Ainverse}), we get
\begin{displaymath}
k(u,x) \bP_{I'} = C^+(u,x) \bP_{I'},
\end{displaymath}
where $C^+(u,x)$ is as in Proposition~\ref{prop:osceledets:GB} (a),
and $I'$ is as in (\ref{eq:def:Iprime}). Similarly, if we may write
\begin{displaymath}
\gamma^-(t) = \bar{k}(v,x) \Lambda^t \bK,
\end{displaymath}
where $\bar{k}(u,x) \in \bK$ and
$\Lambda^t \in \bA_+$. Then, by comparing
(\ref{eq:karlsson:margulis:backwards})  with
(\ref{eq:decomp:A:noinverse}), we get
\begin{displaymath}
\bar{k}(u,x) \bP_{I} = C^-(u,x) \bP_{I'},
\end{displaymath}
where $C^-(u,x)$ is as in Proposition~\ref{prop:osceledets:GB} (b),
and $I$ is as in (\ref{eq:def:I}).
Then, (\ref{eq:opposite:chambers}) implies
(\ref{eq:prop:osceledets:GB:transversal}).
\qed\medskip

\section{Proof of Theorem~\ref{theorem:goodroots}}
\label{sec:proof:theorem:goodroots}

\subsection{Conformal blocks  and Schmidt's criterion} We will use the following criterion of K. Schmidt \cite{Schmidt:Etds1981} for the detection of conformal blocks.

\begin{definition}[cf. Definition 4.6 in \cite{Schmidt:Etds1981}]
\label{def:Schmidt:bounded}
We say that a cocycle $A: G \cross X \to \bH$ is \emph{Schmidt-bounded} if, for every
$\varepsilon>0$, there exists a compact set $\cK(\epsilon) \subset
\bH$ such that
$$\widehat{\beta\times\nu}\left(\left\{((v,u),x)\in\hat{\Omega}\times
    X: A^n(v,u,x) \not\in \cK(\epsilon) \right\}\right)<\varepsilon$$
for all $n\in\mathbb{N}$.
\end{definition}

The importance of this notion in the search of conformal blocks
becomes apparent in view of the next result, which follows from
\cite[Theorem~4.7]{Schmidt:Etds1981}.
\begin{theorem}[Schmidt]
\label{theorem:Schimidt} $A(.,.)$ is
  Schmidt-bounded if and only if there exists a measurable map $C:X\to
  \bH$ such that the cocycle $C(g(x))A(g,x)C(x)^{-1}$ takes its
  values in a compact subgroup of $\bH$.
\end{theorem}

\subsection{Proof of Theorem~\ref{theorem:goodroots}.}
\label{sec:subsec:proof:theorem:goodroots}
We use the notation from \S\ref{sec:subsec:proof:osceledets:a}.
\begin{lemma}
\label{claim:non:singular}
For any $\alpha \in I$, let $\alpha' = - w_0 \alpha w_0^{-1}$ (so that
$\alpha' \in I'$). Then,
$\beta \cross \nu$-almost all $(u,x) \in \Omega \cross X$,  the
measure $\eta_{u,x}^{\alpha'}$ has no atoms; i.e. for any
element $\bar{k}_{u,x} \in \bK$, we have
$\eta_{u,x}^{\alpha'}(\{\bar{k}_{u,x} \bP_{\alpha'} \}) = 0$.
\end{lemma}

\bold{Proof.}
Suppose there exists $\delta > 0$ so that,
for some $\alpha' \in I'$ and for a set $(u,x)$ of positive measure,
there exists $\bar{k}_{u,x} \in \bK$ with
$\eta_{u,x}^{\alpha'}(\{\bar{k}_{u,x} \bP_{\alpha'} \}) > \delta$.
Then this happens for a subset of full measure
by ergodicity. Note that (\ref{eq:measures:converge}) holds.

Then, by Lemma~\ref{lemma:strongly:regular} (b),
for $\beta \cross \nu$ almost all $(u,x)$,
  $\eta^{\alpha'}_{u,x}(\{\bar{k}_{u,x} \bP_{\alpha'} \}) \ge 1-\delta$
  (and thus $\bar{k}_{u,x} \bP_{\alpha'}$ is
  unique) and, as $n \to
\infty$ along $\cN_\delta(u,x)$ (where $\cN_{\delta}(u,x)$ was defined \eqref{eq:Ndelta}), we have:
\begin{displaymath}
\alpha'(a_n(u,x)) \to \infty,
\end{displaymath}
and
\begin{equation}
\label{eq:K:converges}
k_n(u,x) \bP_{\alpha'} \to \bar{k}_{u,x} \bP_{\alpha'},
\end{equation}
Then, by (\ref{eq:d:dbar:conversion}),
\begin{equation}
\label{eq:drds}
\alpha(\bar{a}_n(u,x)) \to \infty,
\end{equation}
and
\begin{displaymath}
\bar{k}_n'(u,x)^{-1} w_0 \, \bP_{\alpha'} \to  \bar{k}_{u,x}
\bP_{\alpha'}.
\end{displaymath}
Therefore for any $\epsilon_1 > 0$ there exists a subset
$H_{\epsilon_1} \subset \Omega
\cross X$ of $\beta \cross \nu$-measure at least $1-\epsilon_1$ such that the convergence in
(\ref{eq:drds}) and (\ref{eq:K:converges})
is uniform over $(u,x) \in H_{\epsilon_1}$.
Hence there exists $M > 0$ such that for any $(u,x) \in H_{\epsilon_1}$,
and any $n \in \cN_\delta(u,x)$ with $n > M$,
\begin{equation}
\label{eq:betaX:small:nbhd:er}
\bar{k}'_n(u,x)^{-1} w_0 \bP_{\alpha'} \in
\Nbhd_{\epsilon_1}(\bar{k}_{u,x} \bP_{\alpha'}).
\end{equation}
By Lemma~\ref{lemma:non:atomic:lyapunov} (see also the proof of
\cite[Lemma~14.4]{Eskin:Mirzakhani:measures})
there exists a subset $H''_{\epsilon_1} \subset X$ with $\nu(H''_{\epsilon_1}) >
1-c_2(\epsilon_1)$ with $c_2(\epsilon_1) \to 0$ as $\epsilon_1 \to 0$ such
that for all $x \in H''_{\epsilon_1}$, and any $g \in \bH$,
\begin{displaymath}
\eta_x( \Nbhd_{2\epsilon_1}(gJ)) < c_3(\epsilon_1),
\end{displaymath}
where $c_3(\epsilon_1) \to 0$ as $\epsilon_1 \to 0$.
Let
\begin{equation}
\label{eq:def:H:prime:epsilon:1}
H'_{\epsilon_1} = \{ (u,x,z) \st (u,x) \in H_{\epsilon_1}, \quad x \in H''_{\epsilon_1} \quad \text{and}
\quad d(z, \bar{k}_{u,x} J) > 2 \epsilon_1 \}.
\end{equation}
Then, $(\beta \cross \hat{\nu})(H'_{\epsilon_1}) > 1-\epsilon_1 - c_2(\epsilon_1) -
c_3(\epsilon_1)$, hence $(\beta \cross \hat{\nu})(H'_{\epsilon_1}) \to
1$ as $\epsilon_1 \to 0$.

We now claim that for $(u,x,z) \in H'_{\epsilon_1}$ and $n \in
\cN_\delta(u,x)$, we have
\begin{equation}
\label{eq:claim:bar:k:n:prime:far}
d(\bar{k}_n'(u,x)z, (\bar{\bN}_\alpha \bP_\alpha)^c) > \epsilon_1.
\end{equation}
Suppose not, then there exist $(u,x,z) \in H'_{\epsilon_1}$ and $n \in
\cN_\delta(u,x)$ such that
\begin{displaymath}
d(\bar{k}_n'(u,x)z, (\bar{\bN}_\alpha \bP_\alpha)^c) \le \epsilon_1.
\end{displaymath}
Let $W_\alpha \subset W$ denote the subgroup of the Weyl group which
fixes $\bM_\alpha$. Then,
\begin{displaymath}
d(\bar{k}_n'(u,x)z, w_0 \bigsqcup_{w \not\in W_\alpha
  w_0^{-1} W_\alpha} \bB w \bB) \le \epsilon_1.
\end{displaymath}
Hence, \mcc{fix inverses?}
\begin{equation}
\label{eq:tmp:d:z:bar:k}
d(z, \bar{k}_n'(u,x)^{-1} w_0 \bigsqcup_{w \not\in
  W_\alpha w_0^{-1} W_\alpha} \bB w \bB) \le \epsilon_1.
\end{equation}
Note that
\begin{displaymath}
\bP_{\alpha'} \bigsqcup_{w \not\in W_\alpha w_0^{-1}
  W_\alpha} \bB w \bB =
\bigsqcup_{w \not\in W_\alpha w_0^{-1} W_\alpha} \bB w \bB.
\end{displaymath}
By (\ref{eq:betaX:small:nbhd:er}) and (\ref{eq:tmp:d:z:bar:k}),
this implies that
\begin{displaymath}
d(z, \bar{k}_{u,x} \bigsqcup_{w \not\in W_\alpha
w_0^{-1} W_\alpha} \bB w \bB)
\le 2\epsilon_1,
\end{displaymath}
contradicting (\ref{eq:def:H:prime:epsilon:1}). This completes the
proof of (\ref{eq:claim:bar:k:n:prime:far}).

Therefore, in view of Lemma~\ref{lemma:sigma:alpha:properties},
there exists $C = C(\epsilon_1)$, such that
for any $(u,x,z) \in H'_{\epsilon_1}$, any $n \in \cN_\delta(u,x)$ with $n > M$,
\begin{equation}
\label{eq:sigma:lower:bound}
\hat{\sigma}_\alpha(A(u_n \dots u_1,x),z) \ge \alpha(A(u_n \dots
u_1,x)).
\end{equation}
By (\ref{eq:drds}) and (\ref{eq:decomp:A:noinverse}), this implies
that for $(u,x,z) \in H'_{\epsilon_1}$,
\begin{equation}
\label{eq:infinite:limit}
\hat{\sigma}_\alpha(A(u_n \dots u_1,x),z) \to \infty \quad\text{as $n
  \to \infty$ along $\cN_\delta(u,x)$.}
\end{equation}
Since $(\beta \cross \hat{\nu})(H'_{\epsilon_1}) \to 1$ as $\epsilon_1
\to 0$, (\ref{eq:infinite:limit}) holds
for $\beta \cross \hat{\nu}$ almost all
$(u,x,z) \in \Omega \cross \hat{X}$.

Let $\sigma_\alpha: \Omega \cross \hat{X} \to \reals$ be defined by
$\sigma_\alpha(u,x,z) = \bar{\sigma}_\alpha(u_1,x,z)$,
where $\bar{\sigma}_{\alpha}$ is as in Lemma~\ref{lemma:furstenberg:formula}.
Then, the left
hand side of (\ref{eq:infinite:limit}) is exactly
\begin{displaymath}
\sum_{j=0}^{n-1} \sigma_\alpha(\hat{T}^j(u,x,z)).
\end{displaymath}
Also, we have $n \in \cN_\delta(u,x)$ if and only if $T^n(u,x) \in \Omega\times
\cK_{\delta}$.
Then, by \cite[Lemma~C.6]{Eskin:Mirzakhani:measures},
\begin{displaymath}
\int_{\Omega \cross \hat{X}} \sigma_\alpha(q)
\, d(\beta \cross \hat{\nu})(q) > 0.
\end{displaymath}
By Lemma~\ref{lemma:furstenberg:formula} (Furstenberg's formula), the
above expression is $\lambda_\alpha$.
Thus $\lambda_\alpha > 0$, contradicting
our assumption that $\alpha \in I$.
This completes the proof of the lemma.
\qed\medskip

\bold{Proof of Theorem~\ref{theorem:goodroots}.}
Let $C^+(u,x) \in \bH$ and $C^-(v,x) \in \bH$ be
as in Proposition~\ref{prop:osceledets:GB}. By
Proposition~\ref{prop:osceledets:GB}(c), for a.e. $(v,u,x)$,
\begin{displaymath}
C^+(u,x)^{-1}C^-(v,x) = p_{I'}(v,u,x) w_0 p_I(v,u,x) \qquad\text{where
  $p_I(v,u,x) \in \bP_I$, $p_{I'}(v,u,x) \in \bP_{I'}$}.
\end{displaymath}
Let
\begin{displaymath}
C_1(v,u,x) = C^+(u,x) p_{I'}(v,u,x) = C^-(v,x)p_I(v,u,x)^{-1} w_0^{-1}.
\end{displaymath}
Then, by Proposition~\ref{prop:osceledets:GB}(a) and (b),
\begin{displaymath}
C_1(\hat{T}^n(v,u,x))^{-1} A^n(v,u,x) C_1(v,u,x) \in
\bP_{I'} \cap w_0 \bP_{I} w_0^{-1} = \bM_{I'} \bA_{I'}.
\end{displaymath}
Let
\begin{displaymath}
A_{I'}^n(v,u,x) := C_1(\hat{T}^n(v,u,x))^{-1} A^n(v,u,x) C_1(v,u,x).
\end{displaymath}
Suppose $\delta > 0$. Then there exist compact sets $\cK_2(\delta)
\subset \hat{\Omega} \cross X$ with $\widehat{\beta \cross
  \nu}(\cK_2(\delta)) > 1-\delta$ and $\cK_3(\delta)\subset \bH$ such
that for $((v,u),x) \in \cK_2(\delta)$, $C_1(v,u,x) \in \cK_3(\delta)$.

Therefore, by (\ref{eq:alpha:in:terms:of:omega:alpha}) and Lemma~\ref{lemma:properties:w:alpha}, there exists $c_1(\delta) \in \reals_+$ such that for
all $((v,u),x) \in \cK_2(\delta)$ and all $n \in \natls$ with
$T^n(v,u,x) \in \cK_2(\delta)$, we have, for all $\alpha \in \Delta$,
\begin{equation}
\label{eq:AIprime:close:to:An}
|\alpha(A_{I'}^n(v,u,x)) - \alpha(A^n(v,u,x))| \le c_1(\delta).
\end{equation}

Let now $\epsilon= \epsilon(\delta)$ and $\cK(\delta)\subset X$ be as in
the proof of
Proposition~\ref{prop:osceledets:GB}(a), so that for $x \in
\cK(\delta)$, the measure $\eta_x$ is $(\epsilon,\delta)$-regular.

By Lemma~\ref{claim:non:singular}, for all $\alpha' \in I'$, the
measures $\eta_{u,x}^{\alpha'}$ are non-atomic. Therefore,
we can find $\epsilon' = \epsilon'(\delta)$ and $\cK'(\delta) \subset
\Omega \cross X$ such that for $(u,x) \in \cK'(\delta)$, and all
$\alpha' \in I'$, for any $z \in \bH/\bP_{\alpha'}$,
\begin{displaymath}
\eta_{u,x}^{\alpha'}(\Nbhd_{\epsilon'}(z)) < \epsilon.
\end{displaymath}
Let
\begin{displaymath}
\cK_1(\delta) = \{(u,x) \in \Omega \cross X \st x \in
\cK(\delta), (u,x) \in \cK'(\delta) \}
\end{displaymath}
Then, by (\ref{eq:measures:converge}), (\ref{eq:decomp:Ainverse}) and
Lemma~\ref{lemma:strongly:regular}(a),
there exists $c_2=c_2(\delta) \in \reals^+$
such that for all $(u,x) \in \cK_1(\delta)$, all  $n$ with $T^n(u,x)
\in \cK_1(\delta)$, and any $\alpha' \in I'$,
\begin{displaymath}
\alpha(A(u_n \dots u_1, x)^{-1}) < c_2(\delta),  \text{where $\alpha =
  - w_0 \alpha' w_0^{-1}$}.
\end{displaymath}
Thus, by (\ref{eq:decomp:A:noinverse}) and
(\ref{eq:d:dbar:conversion}), for all $(u,x) \in \cK_1(\delta)$ and
all  $n \in \natls$ such that $T^n(u,x) \in \cK_1(\delta)$ and all $\alpha' \in I'$,
\begin{equation}
\label{eq:alpha:prime:A:bounded}
\alpha'(A(u_n \dots u_1, x)) < c_2(\delta).
\end{equation}
Let  $\tilde{A}_{I'}^n(v,u,x)$ denote the part of $A_{I'}^n(v,u,x)$
which lies in $\bM_{I'}$. Then,
\begin{equation}
\label{eq:A:Iprime:A:tilde}
\alpha'(\tilde{A}_{I'}^n(v,u,x)) = \alpha'(A_{I'}^n(v,u,x))
\qquad\text{for $\alpha' \in I'$}.
\end{equation}
Let $\cK_1'(\delta) = \cK_2(\delta) \cap \{ (v,u,x) \st (u,x) \in
\cK_1(\delta) \}$.
It follows from (\ref{eq:AIprime:close:to:An}),
(\ref{eq:alpha:prime:A:bounded}) and (\ref{eq:A:Iprime:A:tilde}),
that for all $\alpha'
\in I'$, all $(v,u,x) \in \cK_1'(\delta)$ and all $n \in
\natls$ with $T^n(v,u,x) \in \cK_1'(\delta)$,
\begin{displaymath}
\alpha'(\tilde{A}_{I'}^n(v,u,x)) \le c_3(\delta).
\end{displaymath}
Note that $\widehat{\beta \cross \nu}(\cK_1'(\delta)) > 1-
4\delta$. Since $\delta > 0$ is arbitrary,
it follows that $\tilde{A}_{I'}$ is Schmidt-bounded (see
Definition~\ref{def:Schmidt:bounded}).  Therefore, by
Theorem~\ref{theorem:Schimidt}, there exists $\tilde{C}:
\hat{\Omega} \cross X \to \bM$ such that $\tilde{C}(\hat{T}^n(v,u,x))^{-1}
\tilde{A}_{I'}^n(v,u,x) \tilde{C}(v,u,x) \in \bK \cap
\bM_{I'}$. Let
\begin{displaymath}
C(v,u,x)=C_1(v,u,x) \tilde{C}(v,u,x) w_0.
\end{displaymath}
Then,
\begin{displaymath}
C(T^n(v,u,x))^{-1} A^n(v,u,x) C(v,u,x) \in w_0^{-1} (\bM_{I'} \cap \bK)
\bA_{I'} w_0 = (\bM_{I} \cap \bK) \bA_{I}.
\end{displaymath}
Thus, (\ref{eq:conjugate:into:KM:AI}) holds.

Finally, it is easy to see that (\ref{eq:theorem:goodroots}) follows from
(\ref{eq:AIprime:close:to:An}) and the definition of the
$\lambda_\alpha$ (cf. the argument in the proof of \cite[Lemma 1.5]{Goldsheid:Margulis}).
\qed\medskip

\section{Proof of
Theorem~\ref{theorem:semisimple:twosided:lyapunov}}
\label{sec:proof:theorem:semisimple:twosided:lyapunov}

Let $L$ be a vector space, and suppose $\bH$ is a subgroup of
$SL(L)$. We assume that the action of $\bH$ on $L$ is irreducible,
in the sense that no non-trivial proper subspace of $L$ is fixed by
$\bH$.

Let $\bK$, $I$, $\bA_I$ and $\bM_I$ be as in
Theorem~\ref{theorem:goodroots}. By Theorem~\ref{theorem:goodroots}, we
may assume that the cocycle $A(\cdot, \cdot)$ takes values in $(\bK
\cap \bM_I) \bA_I$. We choose an inner product on $L$ which is
preserved by $\bK$.

Then, the block conformality of
Theorem~\ref{theorem:semisimple:twosided:lyapunov} follows from the
corresponding statement in Theorem~\ref{theorem:goodroots}.


We note that, by Theorem~\ref{theorem:goodroots}, there exists $a_*$ in
the interior of $\bA_I$ such that
the Lyapunov exponents of $A(\cdot, \cdot)$ are given by
expressions of the form $\omega(\log a_*)$, where $\omega$ is a weight of the
action of $\bH$ on $L$.

Let $\cV'_\omega \subset L$ be the  subspace corresponding to the weight
$\omega$; then for $a \in \bA$, and $v \in \cV'_\omega$,
\begin{displaymath}
a \cdot v = \omega(\log a) v.
\end{displaymath}
Let $\omega_0$ be the highest weight. (It exists and has multiplicity
$1$ because the action of $\bH$ on $L$ is irreducible). Then, the top
Lyapunov exponent $\lambda_1$ of $A(\cdot, \cdot)$ is $\omega_0(\log a_*)$.
Then, since the action of $\bH$ on $L$  is
irreducible, the Lyapunov subspace $\cV_1$ of $A(\cdot, \cdot)$
corresponding to the
Lyapunov exponent $\lambda_1$  is given by
\begin{displaymath}
\cV_1 = \bigoplus_{\omega \in S_0} \cV'_\omega,
\end{displaymath}
where $S_0$ consists of weights of the form
\begin{displaymath}
\omega_0 - \sum_{\alpha \in I} c_\alpha \alpha.
\end{displaymath}
Recall that for $\alpha \in I$, $\alpha(\bA_I) = 0$. Therefore, for $a
\in \bA_I$ and $v \in \cV_1$,
\begin{displaymath}
a \cdot v = \omega_0(\log a) v.
\end{displaymath}
Then, for $k \in (\bK \cap \bM_I)$, $a \in \bA_I$,
\begin{equation}
\label{eq:final:thm:almost}
(k a) \cdot v = \omega_0(\log a) k \cdot v.
\end{equation}
Since $A(\cdot, \cdot)$ takes values in $(\bK \cap \bM_I)\bA_I$,
(\ref{eq:top-single-block-twosided}) follows from
(\ref{eq:final:thm:almost}).
\qed

\end{document}